\documentclass[12pt,reqno]{amsart}
\newtheorem{theorem}{Theorem}[section]
\newtheorem{lemma}{Lemma}[section]

\theoremstyle{definition}
\newtheorem{example}[theorem]{Example}

\theoremstyle{remark}
\newtheorem{remark}{Remark}[section]

\numberwithin{equation}{section}
\usepackage{bm,cite}
\usepackage{graphicx}
\usepackage{enumerate}
\usepackage{subfig}
\usepackage{algorithm}  
\usepackage{algorithmic}  
\usepackage{soul}
\begin{document}
\title{Optimal transportation for electrical impedance tomography}


\author{Gang Bao}
\address{School of Mathematical Sciences, Zhejiang University, Hangzhou 310027, China}
\curraddr{}
\email{baog@zju.edu.cn}
\thanks{}
\author{Yixuan Zhang}
\address{School of Mathematical Sciences, Zhejiang University, Hangzhou 310027, China}
\curraddr{}
\email{11935010@zju.edu.cn}
\thanks{}


\date{}

\dedicatory{}

\begin{abstract}
    This work establishes a framework for solving inverse boundary problems with the geodesic based quadratic Wasserstein distance ($W_{2}$). A general form of the Fréchet gradient is systematically derived by optimal transportation (OT) theory. In addition, a fast algorithm based on the new formulation of OT on $\mathbb{S}^{1}$ is developed to solve the corresponding optimal transport problem. The computational complexity of the algorithm is reduced to $O(N)$ from $O(N^{3})$ of the traditional method. Combining with the adjoint-state method, this framework provides a new computational approach for solving the challenging electrical impedance tomography (EIT) problem. Numerical examples are presented to illustrate the effectiveness of our method. 
\end{abstract}

\maketitle
\section{Introduction}
    
The theory of optimal transportation was originally proposed by Monge\cite{monge1781memoire} and later generalized by Kantorovich \cite{kantorovich1960mathematical}. It gives a framework for comparing two probability measures by seeking the minimal cost of rearranging one measure into the other. Optimal transport is closely related to many branches of mathematics, such as partial differential equations, probability analysis, Riemannian geometry, and functional analysis \cite{santambrogio2015optimal}, \cite{villani2021topics}. It has also found applications in a wide range of different fields, including machine learning, economics, optical design, imaging sciences, and graphics \cite{Glimm2003optical}, \cite{Wang2004on}, \cite{haker2004optimal}. Meanwhile, due to the high computational complexity of OT, many numerical algorithms have been developed for practical applications \cite{peyre2019computational}. These algorithms are studied from different perspectives, especially linear programming, the Monge-Ampère equation, and the dynamic formulation of OT. More recently, the Sinkhorn algorithm \cite{Cuturi2013Sinkhorn} has been proposed to solve the entropy regularized OT, which significantly improves computation efficiency.

Over the last few years, optimal transport has been applied to solve inverse problems \cite{abraham2017tomographic}, \cite{chen2018quadratic}, \cite{metivier2016optimal}, \cite{yang2018application}, \cite{heaton2022wasserstein}. A general framework based on OT is to use the Wasserstein distance to measure the discrepancy of datasets in data matching problems. This is an appealing approach since the Wasserstein distance, especially the quadratic Wasserstein distance ($W_{2}$), has the ability to capture both amplitude and spatial information. Compared with the traditional $L^{2}$ metric, $W_{2}$ has better convexity and is more robust to noise \cite{engquist2020quadratic}.  In \cite{engquist2014application}, the $W_{2}$ distance was first introduced to process seismic signals. Subsequently, various types of Wasserstein distance have been applied to the earthquake location problem and the full wave inversion, to mitigate the cycle skipping issues \cite{yang2018application}, \cite{chen2018quadratic}, \cite{metivier2016optimal}, \cite{zhou2018wasserstein}, \cite{metivier2019graph}, \cite{engquist2020optimal}.

This work aims to develop a method based on the quadratic Wasserstein distance to solve severely ill-posed inverse problems. Our particular focus is on the electrical impedance tomography (EIT) problem, also known as Calderón's problem \cite{calderon2006inverse} in the mathematics literature. The problem is to determine the electrical conductivity of a medium from the voltage to current map on the boundary. A typical strategy is to solve the inverse problem with iterative optimization methods, which attempts to minimize certain discrepancy functional \cite{kohn1990numerical}, \cite{cheney1999electrical}, \cite{cheney1990noser}, \cite{wexler1985impedance}. Due to the severe ill-posedness and the nonlinearity of the problem, various regularization strategies have been adopted to resolve the instability \cite{rondi2001enhanced}, \cite{chung2005electrical}, \cite{jin2012reconstruction}. Recently, methods based on deep neural networks have also been applied to solve the EIT problem \cite{bao2020numerical}, \cite{fan2020solving}, especially for the high dimension problems. However, solving the EIT problem in a stable way remains a big challenge in computational inverse problems. Our work here is devoted to the numerical solution of the two-dimensional EIT problem, where the conductivity is located within a disk. It should be pointed out that the two-dimensional EIT problem is particularly challenging mathematically due to the fact that the 2D inverse problem is formally determined. In addition, the two-dimensional EIT problem arises in many practical applications, such as medical imaging \cite{adler2021electrical} and flow monitoring.

For many inverse problems \cite{abraham2017tomographic}, \cite{bao2015inverse}, the observed data are measured on the boundary, which is usually a low-dimensional manifold in the Euclidean space. It is natural to incorporate the geometric information of the manifold into the metric. Instead of using the traditional Euclidean distance as the cost function for OT, in this work, we consider the transportation problem on the manifold and adopt the corresponding geodesic distance as its cost function. This not only improves the computational efficiency, but also better captures the geometric features of the data \cite{solomon2015convolutional}. In particular, the geodesic-based $W_{2}$ distance is employed as the misfit function to solve the two dimensional EIT problem whose data is measured on the circle. Based on our new formulation of OT on $\mathbb{S}^{1}$, an efficient algorithm is designed to calculate the quadratic Wasserstein distance. The complexity of our method is reduced to $O(N)$, while the complexity for the simplex algorithm and the Sinkhorn algorithm is $O(N^{3})$ and $O(N^{2})$, respectively. A crucial step for solving the resulting optimization problem is to develop a new framework for computing the Fréchet gradient of $W_{2}$, which is achieved by observing the explicit connection between the Kantorovich potential and the optimal map. The framework presents a strong contrast to the existing approaches \cite{chen2018quadratic}, \cite{yang2018application}, where the gradient is derived through the perturbation of the fully nonlinear Monge-Ampère equation. Finally, a gradient descent algorithm is employed to solve the optimization problem of EIT.

The paper is organized as follows. In Section 2, by exploring the particular structure and properties of OT under the quadratic cost, we develop a new way to derive the Fréchet gradient of $W_{2}$. Section 3 is devoted to the optimal transportation problem on $\mathbb{S}^{1}$. The simplified formulation of OT on $\mathbb{S}^{1}$ is derived, whose properties are provided to confirm the solvability of this formulation. The corresponding numerical method for solving OT is presented in Section 4. In Section 5, based on our efficient algorithm of computing $W_{2}$ distance and gradient, an adjoint state method is developed to solve the EIT inverse problem. In Section 6, numerical results are provided to demonstrate the effectiveness and efficiency of our method. The paper is concluded with some general remarks in Section 7.

\section{Optimal Transport}
In this section, the prime and dual formulations of optimal transport are presented. Under the quadratic cost, the connection between these two formulations is explored to give rise to a new, straightforward characterization of the $W_{2}$ Fréchet gradient.

Consider two probability measures $\mu$ and $\nu$ defined on complete and separable metric spaces (i.e., polish spaces) $X$ and $Y$, respectively. The cost function $c(x, y)$ maps pairs $(x, y) \in X \times Y$ to $\mathbb{R} \cup\{+\infty\}$. Then Monge's mass transportation problem is to minimize the functional
\begin{equation}\label{OT}
 \int_{X} c(x, T(x))  \mathrm{d} \mu(x)
\end{equation}
over all of the rearrange maps $T$ from $\mu$ to $\nu$:
\begin{equation}\label{MP}
    \Pi(\mu,\nu):=\left\{T: X \rightarrow Y: \forall B \subset Y, \nu(B)=\mu(T^{-1}(B))\right\}.
\end{equation} 
The dual problem of (\ref{OT}) is to maximize
\begin{equation}\label{dual}
 J(\varphi, \psi):=\int_{X} \varphi(x) \mathrm{d} \mu(x)+\int_{Y} \psi(y)\mathrm{d}\nu(y)
\end{equation}
over $\text{Lip}_{c}$, the set of continuous functions $(\varphi,\psi)\in C(X)\times C(Y)$ satisfying
\begin{equation}\label{Lipc}
    \varphi(x)+\psi(y) \leq c(x, y),\quad \forall (x,y)\in X\times Y.
\end{equation}
The standard duality result \cite{villani2021topics} shows that the infimum of (\ref{OT}) is equal to the supremum of (\ref{dual}). The dual formulation is a linear optimization problem under convex constraints, which is desirable for designing numerical algorithms. It also plays an essential role in characterizing the geometrical structure of optimal transportation.

Let $X=Y:=M$. Assume that $d$ is a metric on $M$. The optimal transportation problem naturally defines a distance between probability measures, often referred to as the Wasserstein distance. The p-Wasserstein distance between $\mu$ and $\nu$ is defined by:
\begin{equation}\label{Wp}
    W_{p}(\mu,\nu) = \left(\inf_{T\in \Pi(\mu,\nu)}\int_{M}d(x,T(x))^{p}\mathrm{d}\mu(x) \right)^{\frac{1}{p}},
\end{equation}
which measures the distance between two distributions as the optimal cost of rearranging one distribution into the other. For the rest of the paper, we consider the most common cases for $(M,d)$ in (\ref{Wp}):
\begin{equation}\label{assp}
\begin{aligned}
(\romannumeral1) &\; M = \mathbb{R}^{d} \text{ and } d(x,y) =|x-y| \text{ or }\\
(\romannumeral2) &\; M \text{ is a compact Riemannian manifold, }d \text{ is the geodesic distance on M}.\\
\end{aligned}
\end{equation}
In fact, the focus is primarily on the quadratic Wasserstein distance ($W_{2}$), since the results of $W_{2}$ are most intuitive in both theory and applications. 

To characterize the relationship between the map $T$ in (\ref{OT}) and the dual pair $(\varphi, \psi)$ in (\ref{dual}), it is crucial to analyze further the dual problem of OT. To this end, we introduce the notion of c-transform. For a continuous function $\varphi$ on $M$, its c-transform is defined by
\begin{equation}\label{c-transform}
    \varphi^{c}(y):= \inf _{x \in M} \left\{c(x, y)-\varphi(x)\right\}.
\end{equation}
For any pair $(\varphi, \psi)\in \operatorname{L i p}_{c}$, using the definition of c-transform and the inequality (\ref{Lipc}), we have $\varphi^{c}\geq \psi$ and further $\varphi^{cc}:= (\varphi^{c})^{c}\geq \varphi$. It follows that 
$$J(\varphi^{cc}, \varphi^{c})\geq J(\varphi, \varphi^{c})\geq J(\varphi, \psi).$$
Therefore, the supremum of $J(\varphi, \psi)$ is attained on a smaller set 
\begin{equation}\label{Phi}
   \Phi_{c}:= \left\{(\varphi^{cc}, \varphi^{c}),\varphi\in C(M)\right\}.
\end{equation}
Note that the set $\Phi_{c}$ is well-defined, since it is evident that $(\varphi^{cc})^{c} = \varphi^{c}$. Hence the dual problem (\ref{dual}) only depends on a single variable $\varphi$. 

Assume that there exist an optimal transport map $T\in \Pi(\mu,\nu)$ minimizing (\ref{OT}) and an optimal dual pair $(\varphi, \psi)\in \Phi_{c}$ maximizing (\ref{dual}). Then
\begin{equation}\label{relation}
    \begin{gathered}
        \int_{M}c(x ,T(x))\mathrm{d}\mu(x) = \int_{M}\varphi(x)\mathrm{d}\mu(x) +  \int_{M}\psi(y)\mathrm{d}\nu(y) \\
        = \int_{M}\left(\varphi(x) + \psi(T(x))\right)\mathrm{d}\mu(x).
    \end{gathered}
\end{equation}
This optimal potential $\varphi$ is called the Kantorovich potential. From (\ref{relation}) and (\ref{Lipc}),
\begin{equation}\label{equ1}
    \varphi(x) + \psi(T(x)) = c(x, T(x)),\quad \mathrm{d}\mu\; \mathrm{almost\; everywhere}.
\end{equation}
To illustrate the ideas, consider the Euclidean case $M=\mathbb{R}^{d}$. For the quadratic cost $c(x,y) = \frac{1}{2}d(x,y)^{2}$, after rearranging terms, we obtain from (\ref{Lipc})
$$x\cdot y\leq \left(\frac{1}{2}|x|^{2}-\varphi(x)\right)+\left(\frac{1}{2}|y|^{2}-\psi(y)\right).$$
Thus c-transforms can be converted into Legendre transforms by introducing $\varphi^{*}(x):=\frac{1}{2}|x|^{2}-\varphi(x)$ and $\psi^{*}(y):=\frac{1}{2}|y|^{2}-\psi(y)$. From (\ref{c-transform}) and (\ref{Phi}), the relationship between $\varphi^{*}$ and $\psi^{*}$ is given by the following Legendre transforms:
\begin{equation}\label{legendre}
\begin{gathered}
    \psi^{*}(y):= \sup _{x \in M} \{x\cdot y-\varphi^{*}(x) \}\\
    \varphi^{*}(x):= \sup _{y \in M} \{x\cdot y-\psi^{*}(y)\}
    \end{gathered}
\end{equation}
 Both $\varphi^{*}$ and $\psi^{*}$ are convex because they are defined as the supremum of a family of linear functions. We may assume that $\varphi^{*}$ is differentiable on $M$. The first-order optimality condition of (\ref{legendre}) implies that $ y - \nabla\varphi^{*} (x)= 0$, which motivates us to define a map $t$: 
\begin{equation}
    \left\{\begin{array}{ccc}t:  M & \longrightarrow  &M \qquad\qquad\qquad\\  \quad\; x &\mapsto &\nabla\varphi^{*}(x) = x - \nabla \varphi(x).\end{array}\right.
\end{equation}
Thus for $x\in M$, by the convexity of $\varphi^{*}$, the supremum of (\ref{legendre}) is attained if and only if $y = t(x)$, i.e.
\begin{equation}\label{equ2}
    \varphi(x) + \psi(y) = \frac{1}{2}|x - y|^{2} \quad \mathrm{iff} \quad y = t(x).
\end{equation}
Combining (\ref{equ1}) and (\ref{equ2}), the equation $T(x) = t(x)$ holds $\mathrm{d}\mu$-almost everywhere. Therefore, in this case, the optimal map $T$ has the explicit expression $T(x) = x-\nabla\varphi(x)$ in terms of the Kantorovich potential $\varphi$ under the quadratic cost function. 

This straight-forward derivation provides an important connection between the optimal map and the Kantorovich potential. More generally, the following theorem summarizes the existence and characterization of the optimal transport map:
\begin{theorem}\label{Map}
    (Brenier\cite{brenier1991polar}, McCann\cite{mccann2001polar}) Let $(M,d)$ be the metric space defined in (\ref{assp}) and $c(x,y) = \frac{1}{2}d^{2}(x,y)$ is the quadratic cost. Assume the probability measure $\mu$ is absolutely continuous with respect to the volume measure of $M$. Then there exists a unique solution $T$ to Monge's problem (\ref{OT}), characterized by $T(x) = \exp _{x}(-\nabla \varphi(x))$, where $\varphi:M\rightarrow \mathbb{R}$ is the Kantorovich potential of the dual problem (\ref{dual}), which is unique up to additive constants.
\end{theorem}
In \cite{brenier1991polar}, Brenier gave a rigorous proof of the theorem for the case of the Euclidean space. McCann further generalized the concept to compact Riemannian manifolds in \cite{mccann2001polar}. 

Here in the statement, $\exp$ stands for the exponential map on the tangent bundle $TM$. The notation $\exp_{p}X_{p}$ is the end point of the geodesic segment that starts at $p\in M$ in the direction of $X_{p}\in TM$ with length $|X_{p}|$. In particular, $\exp_{p}X_{p}= x+X_{p}$ in the Euclidean space, corresponding to the result $T(x)=x -\nabla\varphi(x)$. Theorem \ref{Map} assumes that $\mu$ does not give mass to small sets of $M$, which ensures the existence and uniqueness of the optimal map under the quadratic cost. In fact, Monge's optimal transport $T$ may not always exist. A counter-example was given in \cite{santambrogio2015optimal}, where $\mu$ was set to be a weighted sum of Dirac measures. To resolve this issue, Kantorovich proposed an alternative formulation that relaxes the map $T$ to a "multivalued" transport plan \cite{kantorovich1960mathematical}. Here we mainly focus on the Monge problem since we are interested in measures with proper density functions. The corresponding results in this paper can be extended to the Kantorovich problem in a straightforward way. 
\begin{remark}\label{Brenier} (\emph{$c$-cyclical monotonicity})
    In $\mathbb{R}^{d}$, $T(x)= x -\nabla\varphi=\nabla\varphi^{*}$, where $\varphi^{*}$ is a convex function as we discussed before. This result is often referred to as Brenier's theorem \cite{brenier1991polar}. In addition, since $T$ is a gradient of some convex function, $T$ is cyclically monotone \cite{rockafellar1970convex}. That is, for any $\{x_{i}\}_{i=1}^{N} \subset\operatorname{supp}(\mu)$, $\,\sum_{i=1}^{N} x_{i} \cdot T(x_{i}) \geq \sum_{i=1}^{N} x_{i} \cdot T(x_{\sigma(i)})$ for any permutation $\sigma$ on the set $\{1,\cdots, N\}$. The inequality yields
    \begin{equation}\label{cyclical}
        \sum_{i=1}^{N} c\left(x_{i} ,T(x_{i})\right) \leq \sum_{i=1}^{N} c\left(x_{i}, T(x_{\sigma(i)})\right),
    \end{equation}
    where $c$ is the quadratic cost. In fact, (\ref{cyclical}) can be generalized to optimal transportation problems on polish spaces with any continuous cost function $c$. This property of OT is known as "$c$-cyclical monotonicity", which provides alternative arguments for characterizing optimal transport plans; see \cite{ambrosio2013user}, \cite{figalli2011optimal}, \cite{villani2021topics} for details.
\end{remark}

Suppose that measures $\mu$ and $\nu$ have density functions: $\mathrm{d}\mu = f(x)\mathrm{d}x$, $\mathrm{d}\nu = g(y)\mathrm{d}y$, where $\mathrm{d}x$ and $\mathrm{d}y$ are volume elements of $M$. Throughout, the form $W_{p}(f,g)$ instead of $W_{p}(\mu,\nu)$ will be used to indicate the Wasserstein distance between $\mu$ and $\nu$.

\begin{remark}\label{Monge} (\emph{Monge–Ampère equation})
As mentioned in Remark \ref{Brenier}, the optimal transportation map $T(x)=\nabla \varphi^{*}$. Considering the measure-preserving property (\ref{MP}) of $T$, we arrive at the following Monge–Ampère equation using a change of variables technique:
\begin{equation}\label{MA}
    \operatorname{det}\left(D^{2} \varphi^{*}(x)\right)=\frac{f(x)}{g(\nabla \varphi^{*}(x))}. 
\end{equation}
The Caffarelli regularity theorem \cite{caffarelli1996boundary} of (\ref{MA}) shows that if $f, g\in C^{0,\alpha}$ are bounded from above and below by positive constants on their supports and $\operatorname{supp}g$ is convex, then $\varphi^{*}\in C^{2,\alpha}$, i.e., $\varphi\in C^{2,\alpha}$ and $T\in C^{1,\alpha}$.
\end{remark}
Using $W_{2}$ as a misfit function also requires us to access its gradient information. The Fréchet gradient of the Wasserstein distance is related to the corresponding Kantorovich potential by the following theorem:
\begin{theorem}\label{fregradient}
    (Fréchet gradient of the Wasserstein distance) The functional $f \mapsto W_{p}^{p}(f, g)$ is convex, and its subdifferential at $f$ coincides with the set of Kantorovich potentials of (\ref{dual}). If there is a unique Kantorovich potential $\varphi$ up to additive constants, then the Fréchet derivative $\frac{\delta W_{p}^{p}(f,g)}{\delta f} = \varphi$.
\end{theorem}
In fact, Theorem \ref{fregradient} is valid for the minimal transport costs with general continuous cost functions, on top of $c(x,y)=d^{p}(x,y)$. The proof of Theorem \ref{fregradient} is based on a combination of the duality theory and convex analysis; see details in \cite{santambrogio2015optimal}. For $p = 2$, it follows from Theorem \ref{Map} and Theorem \ref{fregradient} that $\frac{\delta W_{p}^{p}(f,g)}{\delta f} = \varphi$. 

For many numerical algorithms\cite{Benamou2000}\cite{benamou2014numerical}, the Kantorovich potential $\varphi$ may not be calculated directly. However, for the quadratic cost, it is fortune that $\varphi$ can be easily obtained using the relation $T = \exp_{x}(\nabla \varphi)$, as described in Remark \ref{gradient_circle}. In fact, it is sufficient to solve either the prime problem for the optimal map $T$ or the dual problem for the Kantorovich potential $\varphi$ to obtain both the value and the gradient of $W_{2}$.

\section{Optimal Transport on $\mathbb{S}^{1}$}
Let $M=\mathbb{R}$. Suppose the probability measures are supported on the interval $[0,1]$. It is well known \cite{villani2021topics} that $W_{2}$ and its optimal transport map are given by
\begin{equation}\label{1d}
W_{2}^{2}(f, g)=\int_{0}^{1}\left|F^{-1}(t)-G^{-1}(t)\right|^{2} \mathrm{~d} t, \quad T(t)=G^{-1}(F(t))
\end{equation}
where $F$ and $G$ are cumulative distributions functions of $f$ and $g$ respectively:
\begin{equation}\label{cumulative}
F(t)=\int_{0}^{t} f(\tau) \mathrm{d} \tau, \quad G(t)=\int_{0}^{t} g(\tau) \mathrm{d} \tau.
\end{equation}
The inverse of the distribution functions are defined by
\begin{equation}\label{inverse-distrib}
F^{-1}(y)=\inf \left\{t: y<F(t)\right\}, \quad G^{-1}(y)=\inf \left\{t: y<G(t)\right\}.
\end{equation}
Formula (\ref{1d}) leads to the algorithm with $O(N)$ complexity. However, for $M=\mathbb{R}^{d}$ with $d\geq 2$, there is no explicit expression available for $W_{2}$. In that case, efficient new approaches are needed to compute the transportation cost since the existing methods directly based on optimization problems (\ref{OT}) and (\ref{dual}) all involve the computational complexity up to $O(N^{3})$. However, as an exception, in the following we show that for $M=\mathbb{S}^{1} \subset \mathbb{R}^{2}$, the optimal transport problem can be reduced to the problem on the real line by cutting the circle at some particular point.

Consider $M=\mathbb{S}^{1} \cong \mathbb{T}=\mathbb{R} /\mathbb{Z}$. For the density function $f$ on $\mathbb{S}^{1}$, we extend its domain from the representative interval $[0,1)$ to $\mathbb{R}$ by making $f(t)$ periodic : $f(t+1)=f(t)$. Thus the definition of distribution function and inverse distribution function are generalized through
\begin{equation}\label{periodic}
    F(t + 1) = F(t) + 1 ,\quad F^{-1}(t+1) =  F^{-1}(t) +1,
\end{equation}
since $f$ and $g$ have unit mass on $[0,1)$. For the cost function $c(x,y)=d^{2}(x,y)$, suppose that $x$, $y$ are the representative elements of $\mathbb{R} /\mathbb{Z}$ belonging to the interval $[0,1)$. Then the geodesic distance $d(x,y)$ on $\mathbb{S}^{1}$ is defined as
$$ d(x,y):=  \operatorname{min}\{|x-y|, 1-|x-y|\}.$$
The exact formula of $W_{2}$ on $\mathbb{S}^{1}$ is stated below:  
\begin{theorem}\label{circle}
Let $f$ and $g$ be two probability distributions on $\mathbb{S}^{1}$, with cumulative distribution functions $F$, $G$ and inverse distribution functions $F^{-1}$, $G^{-1}$ defined by (\ref{cumulative}), (\ref{inverse-distrib}) and (\ref{periodic}), respectively. Let $G^{\alpha}$ denote the function $G + \alpha$. Then the quadratic Wasserstein distance on $\mathbb{S}^{1}$ takes the form:
\begin{equation}\label{W2_circle}
W_{2}^{2}(f, g)=\inf _{\alpha \in \mathbb{R}} \int_{0}^{1}|F^{-1}-\left(G^{\alpha}\right)^{-1}|^{2}\mathrm{d} t.
\end{equation}
Moreover, the optimal map is given by $T(t) = ((G^{\alpha^{*}})^{-1}\circ F)(t)$, where $\alpha^{*}$ is the infimum point in (\ref{W2_circle}).
\end{theorem}
Proof.  We first prove the theorem by assuming $f$ is strictly positive on $\mathbb{S}^{1}$. Let $T$ be the optimal map given in Theorem \ref{Map} and $l(x,y)$ be the geodesic path going from $x$ to $y$, which does not contain $x$ and $y$. We first study the transport patterns of $T$ in order to prove the feasibility of cutting the circle.

Denote $l_{i}:= l(x_{i}, T(x_{i}))$. For any two point $x_{1}$ and $x_{2}$, $x_{1}\neq x_{2}$, one of the following statements must hold:
\begin{itemize}
    \item $l_{1}\cap l_{2} = \emptyset$.
    \item $l_{1}\cap l_{2} \neq \emptyset$, then $l_{1}$ and $l_{2}$ have the same direction, clockwise or both counterclockwise. Furthermore, neither of them is contained in the other.
\end{itemize}
It is sufficient to study the case of $l_{1}\cap l_{2} \neq \emptyset$. Since the cyclical monotonicity (\ref{cyclical}) is valid for any finite sequence on the support of $f$, we have
\begin{equation}\label{assp_c}
    c(x_{1}, T(x_{1}))+ c(x_{2},T(x_{2}))\leq c(x_{1}, T(x_{2}))+ c(x_{2},T(x_{1})).
\end{equation}
Here we prove by contradiction. Assume the directions of $l_{1}$ and $l_{2}$ are different, then $c(x_{1}, T(x_{2}))< c(x_{1},T(x_{1}))$ and $c(x_{2}, T(x_{1}))< c(x_{2},T(x_{2}))$, which contradicts to (\ref{assp_c}). By the convexity of the quadratic cost $c$, similarly we can prove that any path is not contained in other paths. 

Using the statements above, we will show that there is a point $x^{*}$ at which the circle can be cut, that is, there is a point $x^{*}\in \mathbb{S}^{1}$ such that for all $x\neq x^{*} \in \mathbb{S}^{1}$, $x^{*}\notin l(x,T(x))$. 

Again, we prove by contradiction. Assume that for each point $x$, there exists $y\neq x$ such that $x\in l(y,T(y))$. Under this assumption, there is no point $x \in \mathbb{S}^{1}$ such that $x=T(x)$. If $x$ exists, $l(x,T(x))$ is contained in another path, which contradicts the second statement.

Fix a point $x_{0}\in \mathbb{S}^{1}$, $x_{0}\in l(y_{0},T(y_{0}))$ for some $y_{0}\in \mathbb{S}^{1}$. We may assume that $l(y_{0},T(y_{0}))$ goes counter-clockwise. $l_{0}=l(x_{0}, T(x_{0}))$ must move in a counter-clockwise direction since $l_{0}\cap l(y_{0},T(y_{0}))\neq \emptyset$. Denote $x_{1}:= T(x_{0})$. Then $l_{1}=l(x_{1}, T(x_{1}))$ moves in a counter-clockwise direction by a similar argument. Recursively, we obtain a sequence of points $\{x_{i}\}_{i=0}^{\infty}$ with $x_{i+1} = T(x_{i})$. The sequence is strictly increasing in counter-clockwise direction since there is no fixed point. We claim that one of the following situations must hold:
\begin{itemize}
    \item There exists an integer $N$ such that $x_{0}\in l_{N}$.
    \item There exists an integer $N$ such that $x_{N+1} = x_{0}$.
\end{itemize}
If not, the sequence will stagnate before $x_{0}$. Thus it is bounded in the counter-clockwise direction. The limit $x_{\infty}:= \lim_{i\rightarrow \infty}x_{i}$ exists. By the assumption, there is a point $y_{\infty}$ such that $x_{\infty}\in l(y_{\infty}, T(y_{\infty}))$. As $x_{\infty}$ is a limit point, $l_{n}\subset l(y_{\infty}, T(y_{\infty}))$ for large enough $n$, making a contradiction. For
$(x_{0}, x_{1})$, $(x_{1}, x_{2})$, $\cdots$, $(x_{N}, x_{N+1})$, by (\ref{cyclical}),
\begin{equation}\label{assp2}
\sum_{i=0}^{N}c(x_{i}, x_{i+1})\leq \sum_{i=0}^{N}c(x_{i}, x_{\sigma(i+1)}) 
\end{equation}
where $\sigma$ is a permutation of the set $\{1, \cdots, N+1\}$. Specifically, we set 
$$
\sigma(i+1)= 
\begin{cases}N +1, & i=0 \\ 
    i, & 1 \leq i \leq N,
\end{cases} 
$$
For the first situation, $c(x_{0},x_{N+1})< c(x_{N},x_{N+1})$ since $x_{0}$ is in $l_{N}$. For the second situation, $c(x_{0},x_{N+1}) =0$. Hence
$\sum_{i=0}^{N}c(x_{i}, x_{\sigma(i+1)}) = c(x_{0},x_{N+1}) < \sum_{i=0}^{N}c(x_{i}, x_{i+1}), $
which contradicts to (\ref{assp2}). The existence of $x^{*}$ is established.

We can thus cut $\mathbb{S}^{1}$ at $x^{*}$ and reduce the transport problem on the circle to the transport problem on the real line, since all the geodesic paths of the optimal transport map are the same side of $x^{*}$. Taking $x^{*}$ as the new reference point, by (\ref{1d}) and the fact that $(F-c)^{-1}(t) = F^{-1}(t+c)$ for any constant $c$, the optimal transport cost on line $[x^{*},x^{*}+1)$ is 
$$
\begin{gathered}
\int_{x^{*}}^{x^{*}+1}\left|F_{x^{*}}^{-1}(t)-G_{x^{*}}^{-1}(t)\right|^{2}\mathrm{d}t = \int_{x^{*}}^{x^{*}+1}\left|F^{-1}(t+F(x^{*}))-G^{-1}(t+G(x^{*}))\right|^{2} \mathrm{d} t\\
=\int_{x^{*} + F(x^{*})}^{x^{*} + F(x^{*})+1}\left|F^{-1}(t)-G^{-1}(t-\alpha)\right|^{2}\mathrm{d}t
 = \int_{0}^{1}\left|F^{-1}(t)- (G^{\alpha})^{-1}(t)\right|^{2}\mathrm{d}t
\end{gathered}
$$ 
where $F_{x^{*}}(t) = F(t) - F(x^{*})$ , $G_{x^{*}}(t) = G(t) - G(x^{*})$ and $\alpha = F(x^{*}) - G(x^{*})$. The last step follows from the periodic property (\ref{periodic}). 

Calculating the optimal cost of all possible cuttings, we arrive at
$$
W_{2}^{2}(f,g) \geq  \inf_{x^{*}\in\mathbb{S}^{1}}\int_{x^{*}}^{x^{*}+1}\left|F_{x^{*}}^{-1}(t)-G_{x^{*}}^{-1}(t)\right|^{2}\mathrm{d}t \geq \inf_{\alpha\in\mathbb{R}}\int_{0}^{1}\left|F^{-1} - (G^{\alpha})^{-1}\right|^{2}\mathrm{d}t.$$
Constructing the map $T(t) = (G^{\alpha})^{-1}\circ F(t)$ on the circle, we can get 
$$ W_{2}^{2}(f,g) \leq \int_{0}^{1}c\left(t,T(t)\right)f(t)\mathrm{d}t \leq \int_{0}^{1}\left|t-T(t)\right|^{2}f(t)\mathrm{d}t =\int_{0}^{1}\left|F^{-1} - (G^{\alpha})^{-1}\right|^{2}\mathrm{d}t.$$
Taking the infimum with respect to $\alpha$, the identity (\ref{W2_circle}) is proved for $f> 0$. 

For $f\geq 0$ on $\mathbb{S}^{1}$, we have positive density sequences $\{f_{n}\}$ such that $f_{n}\rightarrow f$ pointwise. By the weak convergence of the Wasserstein distance \cite{villani2021topics}, $W_{2}^{2}(f_{n},g)\rightarrow W_{2}^{2}(f,g)$. Regarding the formula on the right side of (\ref{W2_circle}), the convergence can also be obtained from the fact that the optimal $\alpha^{*}$ depends continuously on $f$. 
$\square$\\

We remark that under a more general setting, (\ref{W2_circle}) has also been derived by a different approach \cite{delon2010fast}, namely Aubry–Mather theorem. Our proof here from a different point of view is based on the c-cyclical monotonicity, which is more direct. The idea is inspired by the technique of optimal permutation problem on the circle \cite{rabin2011transportation}. 

Consequently, when $M=\mathbb{S}^{1}$, the Wasserstein distance (\ref{Wp}) can be simplified to (\ref{W2_circle}). To solve (\ref{W2_circle}) effectively, it is important to make the following assumption:
\begin{equation}\label{pset}
    D:=\left\{f\in L^{1}[0,1]: \int_{0}^{1}f\mathrm{d}x =1, \,f\geq \eta \text{ on } \mathbb{S}^{1} \text{ for some } \eta >0 \right\}.
\end{equation}
That is, the density functions are assumed to have a positive lower bound. This assumption is reasonable for our applications. In fact, under the assumption, for $f$, $g\in D$, distribution functions $F$ and $G$ are strictly increasing, and thus their inverse functions (\ref{inverse-distrib}) exist in the classical sense. In addition, it is also required for Caffarelli's regularity theorem of OT discussed in Remark \ref{Monge}.

The following lemma provides an intuitive way to find the infimum point $\alpha^{*}$ in (\ref{W2_circle}).
\begin{lemma}\label{property}
Let $f, g\in D$ be fixed probability density functions on $\mathbb{S}^{1}$. Define 
$$ I(\alpha):=I(\alpha; f,g) := \int_{0}^{1}|F^{-1}(t)-\left(G^{\alpha}\right)^{-1}(t)|^{2}\mathrm{d} t .$$
Then the following properties hold for $I(\alpha)$:
\begin{enumerate}[(i)]
    \item $I(\alpha)$ is a strictly convex function about $\alpha$. 
    \item $I(\alpha)$ is twice differentiable with respect to $\alpha$ on $\mathbb{R}$ and
    \begin{equation}\label{first_der}
        I^{\prime}(\alpha) = 1- 2\int_{0}^{1}F^{-1}(G(t)+ \alpha )\,\mathrm{d} t .
    \end{equation}
    \begin{equation}\label{second_der}
        I^{\prime\prime}(\alpha)=\int_{0}^{1}-\frac{2}{f(F^{-1}(G(t)+\alpha))}\mathrm{d}t = \int_{0}^{1}-\frac{2}{g(G^{-1}(F(t)-\alpha))}\mathrm{d}t.
    \end{equation}
    \item The global minimum of $I(\alpha)$ is uniquely attained on the interval $(-1, 1)$.
\end{enumerate}
\end{lemma}
Proof. Denote $c_{f}^{g}(x,y): = |F^{-1}(x) - G^{-1}(y) |^{2}$, then $I(\alpha; f, g)= \int_{0}^{1}c_{f}^{g}(t, t -\alpha)\mathrm{d} t$. 

(\romannumeral1) The proof of the convexity follows from the idea in \cite{delon2010fast}. Let $\alpha_{1}<\alpha_{2}$, denote $\alpha=\frac{1}{2}\left(\alpha_{1}+\alpha_{2}\right)$. Making the change of variables $t^{{\prime}} = t +\alpha - \alpha_{2}$ and taking into account the periodic structure of $F^{-1}$ and $G^{-1}$, we have
$$
\begin{gathered}
    I(\alpha_{1})= \int_{0}^{1}c_{f}^{g}(t, t -\alpha_{1})\mathrm{d} t, \quad I(\alpha_{2})= \int_{0}^{1}c_{f}^{g}(t, t -\alpha_{2})\mathrm{d} t = \\\int_{\alpha_{2}- \alpha}^{1+ \alpha_{2}- \alpha}c_{f}^{g}(t^{\prime}+ \alpha_{2}- \alpha, t^{\prime} -\alpha)\mathrm{d} t^{\prime} 
    = \int_{0}^{1}c_{f}^{g}(t^{\prime}+ \alpha_{2}- \alpha, t^{\prime} -\alpha)\mathrm{d} t^{\prime},\\
    I(\alpha)=\int_{0}^{1} c_{f}^{g}(t, t -\alpha)\mathrm{d} t = \int_{0}^{1} c_{f}^{g}(t^{\prime}+ \alpha_{2}- \alpha, t^{\prime} -\alpha_{1})\mathrm{d} t^{\prime}.
\end{gathered}
$$
Note that for $c(x,y) = |x-y|^{2}$, it satisfies $c\left(x_{1}, y_{1}\right)+c\left(x_{2}, y_{2}\right)<c\left(x_{1}, y_{2}\right)+c\left(x_{2}, y_{1}\right)$ for all $x_{1}<x_{2}$ and $y_{1}<y_{2}$ due to the convexity of the quadratic distance. Since $F^{-1}(t)< F^{-1}(t + \alpha_{2}- \alpha)$ and $G^{-1}(t- \alpha)<G^{-1}(t-\alpha_{1})$,
$$
\begin{gathered}
2I(\alpha) - I(\alpha_{1}) - I(\alpha_{2}) = \int_{0}^{1} c_{f}^{g}(t, t -\alpha) + c_{f}^{g}(t+ \alpha_{2}- \alpha, t -\alpha_{1})\\
-  c_{f}^{g}(t, t -\alpha_{1}) - c_{f}^{g}(t + \alpha_{2}- \alpha, t -\alpha)\mathrm{d} t<0.
\end{gathered}
$$
Hence $I(\frac{1}{2}(\alpha_{1} + \alpha_{2}))< \frac{1}{2}(I(\alpha_{1})+ I(\alpha_{2}))$. $I(\alpha)$ is strictly convex.

(\romannumeral2) 
Both $F^{-1}$ and $G^{-1}$ are Lipschitz continuous due to the strict positivity of $f$ and $g$. Thus $c_{f}^{g}(x,y)$ is Lipschitz continuous in both variables on any bounded domain, and its partial derivatives exist almost everywhere. As a consequence, 
$$D(t;\alpha, \delta \alpha):= \frac{c_{f}^{g}(t, t - \alpha - \delta \alpha) - c_{f}^{g}(t, t -\alpha)}{\delta \alpha},\quad t\in [0,1]$$
is uniformly bounded for any small perturbation $\delta\alpha$ when $\alpha\in \mathbb{R}$ is fixed. By the dominated convergence theorem,
$$
\begin{gathered}
\lim_{\delta \alpha\rightarrow 0}\frac{I(\alpha + \delta \alpha) - I(\alpha)}{\delta \alpha}=\int_{0}^{1}\lim_{\delta \alpha \rightarrow 0} D(t;\alpha, \delta \alpha)\mathrm{d} t \\
=\int_{0}^{1}-\frac{\partial c_{f}^{g}}{\partial y}(t,t - \alpha )\mathrm{d} t = \int_{0}^{1}-\frac{2}{g(G^{-1}(t-\alpha))}\left(F^{-1}(t) - G^{-1}(t-\alpha)\right)\mathrm{d}t\\
=2\int_{G^{-1}(-\alpha)}^{1 + G^{-1}(-\alpha)}\left(y - F^{-1}(G(y)+ \alpha)\right)\mathrm{d} y= 2\int_{0}^{1}\left(y - F^{-1}(G(y)+ \alpha) \right)\mathrm{d} y\\
=1 - 2\int_{0}^{1}F^{-1}(G(t)+ \alpha)\mathrm{d} t.
\end{gathered}
$$
As for the second derivative, by a similar argument, we obtain
$$
\begin{gathered}
I^{\prime\prime}(\alpha)= \lim_{\delta \alpha\rightarrow 0}\frac{I^{\prime}(\alpha + \delta \alpha) - I^{\prime}(\alpha)}{\delta \alpha} =\int_{0}^{1}-\frac{2}{f(F^{-1}(G(t)+\alpha))}\mathrm{d}t \\
=\int_{0}^{1}-\frac{2}{g(G^{-1}(y))f(F^{-1}(y+\alpha))}\mathrm{d}y=\int_{0}^{1}-\frac{2}{g(G^{-1}(F(t)-\alpha))}\mathrm{d}t.
\end{gathered}    
$$
(\romannumeral3) It is easy to check that 
$$
I^{\prime}(1)=1-2\int_{0}^{1}F^{-1}(G(t)+1)\mathrm{d}t< 1- 2F^{-1}(G(0)+1)<0
$$
and $I^{\prime}(-1)>0$. Then (\romannumeral1) and (\romannumeral2) together imply that there exists a unique $\alpha^{*}\in (-1, 1)$ such that $I^{\prime}(\alpha^{*}) = 0$. It follows from the convexity that $\alpha^{*}$ is the global minimum point. $\square$
\begin{remark}(\emph{Fréchet gradient of $W_{2}$ on $\mathbb{S}^{1}$})\label{gradient_circle}
    According to Theorem \ref{fregradient} , $\frac{\delta W_{2}^{2}(f,g)}{\delta f}$ is the Kantorovich potential of the dual problem. Combining with Theorem \ref{Map}, the Fréchet gradient of $W_{2}^{2}$ on $\mathbb{S}^{1}$ is given by the following integral along the circle:
    \begin{equation}\label{fre_exp}
        \frac{\delta W_{2}^{2}(f,g)}{\delta f} =2 \int_{0}^{t} \left(\tau -  T(\tau)\right)\mathrm{d}\tau + c,
    \end{equation}
    where $T(\tau) = G^{-1}(F(\tau)- \alpha^{*})$ is the optimal transport map between $f$ and $g$ , and $c$ is an arbitrary constant. 
\end{remark}
Lemma \ref{property} indicates that the computation of $W_{2}$ is equivalent to solving the nonlinear equation $I'(\alpha)=0$. The corresponding numerical algorithm is proposed in the next section.

\section{Numerical Method}
For the computation purpose, the density function $f$ defined on $[0,1)$ is discretized on the nodes $\tau_{i} = i*h$, $i=0,1,2,\cdots,N-1$ where $h = \frac{1}{N}$. For convenience, we extend the interval $[0,1)$ to $[-1, 2]$, which is discretized by 
$-1 = t_{-N}<t_{-(N-1)}<\cdots t_{0}<\cdots< t_{2N}< t_{2N +1} =2,$
with $t_{i} = (2 i- 1)*\frac{h}{2}$, $i = -N<i<2N+1$. 

Define on the interval $[-1,2)$:
 $$f_{h}(t) = f_{i}:= \frac{1}{m}f(\tau_{j(i)})\;\, \mathrm{for} \;t\in I_{i},\;\,j(i) = i\,\mathrm{mod}\,N, \quad i =-N, \cdots, 2N$$ 
where $I_{i} := [t_{i}, t_{i+1})$ and the rescale parameter $m=\sum_{i=0}^{N-1}f(\tau_{i})h$ is the mass over the period. Apparently, $f_{h}(t)$ is a periodic piecewise constant function with unit mass in each period. According to (\ref{cumulative}), the cumulative distribution function and inverse cumulative distribution function are given by 
\begin{equation}\label{dis_culmulative}
\begin{gathered}
F_{h}(t) = f_{i}(t - t_{i})+ F_{i}\, ,\;\;t\in I_{i},\qquad F_{h}^{-1}(y) = \frac{y- F_{i}}{f_{i}} + t_{i}\,,\;\; y \in [F_{i}, F_{i+1})
\end{gathered}
\end{equation}
where $F_{-N} = -1$, $F_{-(N-1)}=-1+\frac{h}{2}f_{-N}$, and $F_{i} = F_{i-1} + h*f_{i-1}$ for $-(N-1)<i\leq 2N$. 

The density function $g$ can be discretized in same way to get $g_{h}$ and $G_{h}$.

We are now ready to compute the discretized version of (\ref{first_der}) to solve the optimal transportation problem. The difficulty lies in the fact that (\ref{first_der}) involves the inverse of $F$ in a composite form with $G$, which requires the correspondance of the nodes for $F^{-1}$ and $G$. For completeness, we provide the details below.

Given $\alpha\in(-1,1)$, there exists an integer $i_{\alpha}$ such that $\alpha \in [F_{i_{\alpha}}, F_{i_{\alpha}+1})$, hence $\alpha + 1\in [F_{i_{\alpha} + N}, F_{i_{\alpha} + N+1})$. For two increasing sequences $\{H_{i}^{0}:= G_{i}\}_{i=1}^{N}$ and $\{H_{i}^{1}:=F_{i_{\alpha} +i}-\alpha\}_{i=1}^{N}$,  we sort their values into one increasing sequence, which is denoted by $\{H_{n}\}_{n=1}^{2N}$. The sorting process automatically defines a bijective map $\sigma: (i,j)\rightarrow n$ such that $H_{i}^{j}$ is reordered as $H_{n}$ in the new sequence. Now we define the indexing sequences $\{l_{n}^{0}\}_{n = 1}^{2N}$ and $\{l_{n}^{1}\}_{n = 1}^{2N}$ by:
$$
    \left\{\begin{array}{l}
    l_{n}^{j} =  i \quad  \mathrm{where} \; (i,j) =  \sigma^{-1}(n)\\
    l_{n}^{1-j} = n - l_{n}^{j}
    \end{array}\right.
$$
In consequence, the sequence of nodes $\{T_{n}:=G_{h}^{-1}(H_{n})\}_{n=1}^{2N}$ is easy to compute:
$$T_{n}=\left\{\begin{array}{cl}
    \frac{H_{n}-G_{l_{n}^{0}}}{g_{l_{n}^{0}}} + t_{l_{n}^{0}} & ,\;\text { if } j=1 \\
    t_{i} & ,\;\text { if } j=0
    \end{array} \quad\mathrm{where} \; (i, j)=\sigma^{-1}(n)\right.$$
Adding $T_{0} =0 $ and $T_{2N+1} = 1$ to the sequence, we can get $\{T_{n}\}_{n=0}^{2N+1}$. Reset $\{l_{n}^{1}:=l_{n}^{1}+i_{\alpha}\}_{n=1}^{2N}$. Then for $t\in [T_{n}, T_{n+1})$, 
$$F_{h}^{-1}(G_{h}(t)+\alpha) = K_{n}t + B_{n}, $$
which is a piecewise linear function on $[0,1)$. The parameters $K_{n}$ and $B_{n}$ are computed as
$$K_{n} =\frac{g_{l_{n}^{0}}}{f_{l_{n}^{1}}}, \quad B_{n} =\frac{\alpha +G_{l_{n}^{0}}-g_{l_{n}^{0}}t_{l_{n}^{0}}-F_{l_{n}^{1}}}{f_{l_{n}^{1}}} + t_{l_{n}^{1}}.$$
Finally, the integral (\ref{first_der}) is discretized as  
\begin{equation}\label{dis_exp}
    I_{h}^{\prime}:= 1-2\int_{0}^{1}F_{h}^{-1}(G_{h}(t)+\alpha) \mathrm{d}t = 1-2\sum_{n=0}^{2N}\left(\frac{1}{2}K_{n}(T_{n+1}^{2} - T_{n}^{2})+ B_{n}(T_{n+1}-T_{n})\right).
\end{equation}
The second derivative $I_{h}^{\prime\prime}$ can be computed in the same fashion.
\begin{algorithm}
    \caption{Newton's method}
    \label{alg1} 
    \begin{algorithmic}  
    \REQUIRE $f_{h}$, $g_{h}$ and the precision $\epsilon$.  
     \STATE Initially set $n=0$, $\alpha_{0}=1$. Compute parameters in $F_{h}$ and $G_{h}$.
     \WHILE{$|\alpha_{n}-\alpha_{n-1}|\geq\epsilon$ or $n=0$}
     \STATE Compute $I_{h}^{\prime}(\alpha_{n})$ and $I_{h}^{\prime\prime}(\alpha_{n})$.
     \STATE Update $\alpha_{n+1} :=\alpha_{n}-\frac{I_{h}^{\prime}(\alpha_{n})}{I_{h}^{\prime\prime}(\alpha_{n})}$.
     \STATE \qquad \qquad $n:=n+1$.
     \ENDWHILE
     \STATE Using the output value $\alpha_{h}$, the quadratic Wasserstein distance $W_{2}^{2}(f_{h}, g_{h})$ and the Fréchet gradient can be computed in the same way as the calculation of (\ref{dis_exp}).
    \end{algorithmic}   
\end{algorithm}
\begin{remark}
    Here, Newton's method is used to solve the nonlinear equation $I_{h}^{\prime}=0$. A detailed description of the algorithm is given in Algorithm \ref{alg1}. The strict convexity of $I(\alpha)$ guarantees the algorithm converge to the minimum point. Newton's method takes at most $O(\log\log(\frac{1}{\epsilon}))$ steps to obtain $\alpha$ within accuracy $\epsilon$. Since sorting two increasing sequences requires $2N$ comparisons, each step of evaluating $I_{h}^{\prime}$ and $I_{h}^{\prime\prime}$ takes at most $O(N)$ operations. Consequently, the computational complexity of this algorithm is $O(N\log\log(\frac{1}{\epsilon}))$. However, if we compute $W_{2}$ directly in $\mathbb{R}^{2}$, the computation cost would be extremely high. 
\end{remark}
From Lemma \ref{property}, for every pair $(f,g)\in D$, there is a unique $\alpha\in(-1,1)$ such that $I^{\prime}(\alpha;f,g)=0$. Thus $\alpha$ can be viewed as a function $\alpha(f,g)$.
Discretizations of the density functions lead to errors in $\alpha$ along with $W_{2}$. Next we provide a stability estimate of $\alpha(f,g)$. 

\begin{lemma}\label{estimate}
    (Stability Estimate) 
    Suppose that $\alpha(f,g)$ is the implicit function defined by $I^{\prime}(\alpha;f,g)=0$, and $f_{i}, g_{i} \in D$ are continuous differentiable, $i=1,2$. Then for $\alpha_{i}:= \alpha(f_{i},g_{i})$, the following estimate holds:
    \begin{equation}\label{estimate2}
        |\alpha_{1} - \alpha_{2}|\leq \frac{1}{2}\left(\|f_{1} - f_{2}\|_{L^{1}[0,1]}+\|g_{1} - g_{2}\|_{L^{1}[0,1]}\right).
    \end{equation}
\end{lemma}

Proof. Let $L(\alpha,f,g):= \int_{0}^{1}F^{-1}(G(t)+\alpha)\mathrm{d}t$, then $I^{\prime}(\alpha; f, g) =1- 2L(\alpha, f,g)$. For $(\alpha, f, g)$ with $\alpha = \alpha(f,g)$, we perturb $f$ by an amount $\delta f$ and investigate the resulting change in $L(\alpha, f, g)$ as a functional of $f$. For simplicity, $\tilde{F}(t)$ is used to denote the function $(F + \delta F)(t):=\int_{0}^{t}(f + \delta f)(\tau)\mathrm{d}\tau$. Let $y= F(t)$ and thus $t = F^{-1}(y)$. Applying the Taylor expansion of $\tilde{F}^{-1}$ at $\tilde{F}(t)$, we arrive at
\begin{equation}\label{e_inverse}
\begin{gathered}
    \tilde{F}^{-1}(y) - F^{-1}(y) =  \tilde{F}^{-1}(F(t)) -  \tilde{F}^{-1}(\tilde{F}(t))\\= -\frac{1}{(f+\delta f)(t)}\int_{0}^{t}\delta f(\tau) \mathrm{d}\tau + O((\delta F)^{2})\\=
    -\frac{\int_{0}^{t}\delta f(\tau) \mathrm{d}\tau}{f(t)} + O((\delta F)^{2}),
\end{gathered}
\end{equation}
which holds for almost all $y\in\mathbb{R}$. Omitting high-order terms, we have
$$
\begin{gathered}
    L(\alpha, f + \delta f,g) - L(\alpha, f, g) = \int_{0}^{1}\tilde{F}^{-1}(G(t)+\alpha) - F^{-1}(G(t)+\alpha)\mathrm{d}t\\
    = \int_{0}^{1}\frac{1}{g(G^{-1}(y-\alpha))}\left(\tilde{F}^{-1}(y) - F^{-1}(y)\right)\mathrm{d}y\\
    = \int_{0}^{1}\frac{-\int_{0}^{t}\delta f(\tau) \mathrm{d}\tau}{g(G^{-1}(F(t)-\alpha))}\mathrm{d}t = \int_{0}^{1}\left(\int_{t}^{1}\frac{-1}{g\circ T(\tau)}\mathrm{d}\tau\right)\delta f(t)\mathrm{d}t,
\end{gathered}
$$
where $T=(G^{\alpha})^{-1}\circ F$ is the optimal transport map between $f$ and $g$. Since $\int_{0}^{1}\delta f(t)\mathrm{d}t =0$, the partial derivative $\frac{\partial L}{\partial f}$ at $(\alpha, f,g )$ is 
$$\frac{\partial L}{\partial f} = \int_{0}^{t}\frac{1}{g\circ T(\tau)}\mathrm{d}\tau + c.$$
Similarly, we compute the partial derivative $\frac{\partial L}{\partial g}$ at $(\alpha, f, g)$:
$$\frac{\partial L}{\partial g} =\int_{0}^{1}\frac{\int_{0}^{t}\delta g(\tau)\mathrm{d}\tau}{f(F^{-1}(G(t)+\alpha))}\mathrm{d}t= -\int_{0}^{t}\frac{1}{f\circ T^{-1}(\tau)}\mathrm{d}\tau + c.$$
The partial derivative $\frac{\partial L}{\partial \alpha}$ at $(\alpha, f, g)$ is obtained directly from (\ref{second_der}) in Lemma \ref{property}:
$$\frac{\partial L}{\partial \alpha} = \int_{0}^{1}\frac{1}{g\circ T(t)}\mathrm{d}t = \int_{0}^{1}\frac{1}{f\circ T^{-1}(t)}\mathrm{d}t.$$
Using the chain rule for the equation $L(\alpha,f,g) =\frac{1}{2}$, we obtain
$$\frac{\partial \alpha}{\partial f}\big|_{(f,g)} = -(\frac{\partial L}{\partial \alpha})^{-1}\frac{\partial L}{\partial f} =-\frac{\int_{0}^{t}\frac{1}{g\circ T(\tau)}\mathrm{d}\tau+c}{\int_{0}^{1}\frac{1}{g\circ T(\tau)}\mathrm{d}\tau},\; \frac{\partial \alpha}{\partial g}\big|_{(f,g)} = -(\frac{\partial L}{\partial \alpha})^{-1}\frac{\partial L}{\partial g} = \frac{\int_{0}^{t}\frac{1}{f\circ T^{-1}(\tau)}\mathrm{d}\tau+c}{\int_{0}^{1}\frac{1}{f\circ T^{-1}(\tau)}\mathrm{d}\tau}.$$
For the estimate about $\frac{\partial \alpha}{\partial f}$,
$$\inf_{c\in \mathbb{R}}\left\|\frac{\partial \alpha}{\partial f}\big|_{(f,g)}(t;c)\right\|_{L^{\infty}[0,1]}= \left\| -\frac{\int_{0}^{t}\frac{1}{g\circ T(\tau)}\mathrm{d}\tau}{\int_{0}^{1}\frac{1}{g\circ T(\tau)}\mathrm{d}\tau}+\frac{1}{2}\right\|_{L^{\infty}[0,1]}=\frac{1}{2},$$
where the infimum is achieved at $c^{*} =-\frac{1}{2}\int_{0}^{1}\frac{1}{g\circ T(t)}\mathrm{d}t$. The same estimate is obtained for $\frac{\partial \alpha}{\partial g}$ as well. By the mean value theorem and Hölder's inequality, 
$$
\begin{gathered}
    |\alpha_{1} -\alpha_{2}| = |\alpha(f_{1},g_{1}) - \alpha(f_{2}, g_{2})|\leq |\alpha(f_{1},g_{1}) - \alpha(f_{2}, g_{1})|+ |\alpha(f_{2}, g_{1}) - \alpha(f_{2}, g_{2}) | \leq\\
     \sup\limits_{s\in[0,1]}\left|\int_{0}^{1}(f_{1}-f_{2})\cdot\frac{\partial \alpha}{\partial f}\big|_{(f_{1}+ s(f_{2}-f_{1}),g_{1})}\mathrm{d}t\right| + 
     \sup\limits_{s\in[0,1]}\left|\int_{0}^{1}(g_{1}-g_{2})\cdot\frac{\partial \alpha}{\partial g}\big|_{((f_{2}, g_{1}+s(g_{2}-g_{1}))}\mathrm{d}t\right| \leq\\
     \inf\limits_{c\in \mathbb{R}}\left\|\frac{\partial \alpha}{\partial f}\big|_{(f_{1}+ s(f_{2}-f_{1}),g_{1})}\right\|_{L^{\infty}[0,1]}\|f_{1}-f_{2}\|_{L^{1}[0,1]} + \inf\limits_{c\in \mathbb{R}}\left\|\frac{\partial \alpha}{\partial g}\big|_{((f_{2}, g_{1}+s(g_{2}-g_{1}))}\right\|_{L^{\infty}[0,1]}\|g_{1}-g_{2}\|_{L^{1}[0,1]}\\
    = \frac{1}{2}\|f_{1} - f_{2}\|_{L^{1}[0,1]}+\frac{1}{2}\|g_{1} - g_{2}\|_{L^{1}[0,1]}.   \quad\square
\end{gathered}
$$

\section{Electrical Impedance Tomography}
Let $\Omega$ be an open-bounded domain in $\mathbb{R}^{2}$ with a smooth boundary $\partial \Omega$ and $\sigma\in L^{\infty}(\Omega)$ be strictly positive on $\Omega$. In our problem, $\Omega$ is the unit disk. The EIT forward problem is modeled by the elliptic partial differential equation:
\begin{equation}\label{EIT_Forward}
\begin{aligned}
-\nabla \cdot(\sigma \nabla u) &=0 \quad \text { in } \Omega \\
\sigma \frac{\partial u}{\partial n} &=j \quad \text { on } \partial \Omega
\end{aligned}
\end{equation}
where $u$ and $j$ denote the electrical potential and current, respectively. For each $j\in \tilde{H}^{-\frac{1}{2}}(\partial \Omega):=\{v \in H^{-\frac{1}{2}}(\partial \Omega): \int_{\partial \Omega} v \mathrm{~d} s=0\}$, there is a unique $u\in \tilde{H}^{1}(\Omega):=\{v \in H^{1}(\Omega): \int_{\partial \Omega} v \mathrm{~d} s=0\}$ solves the equation (\ref{EIT_Forward}).
Therefore, for each $\sigma$ satisfying the condition, define a Neumann to Dirichlet operator $\Lambda(\sigma): \tilde{H}^{-\frac{1}{2}}(\partial \Omega)\rightarrow \tilde{H}^{\frac{1}{2}}(\partial \Omega)$:
\begin{equation}
    \Lambda({\sigma})j = \phi,\quad \mathrm{where}\;\phi = \gamma u , 
\end{equation}
where $\gamma$ is the trace operator projecting $\tilde{H}^{1}(\Omega)$ to $\tilde{H}^{\frac{1}{2}}(\partial \Omega)$. The boundary operator $\Lambda({\sigma})$ is also known as NtD map. This map is self-adjoint and positive definite. While the forward problem is calculating $\Lambda(\sigma)$ with known $\sigma$, the inverse problem (EIT) is to reconstruct $\sigma$ from the knowledge of $\Lambda(\sigma)$, which can be formulated as an optimization problem:
\begin{equation}\label{optimization}
    \begin{gathered}
    \sigma^{*} =  \mathop{\mathrm{argmin}}\limits_{\sigma \in \mathcal{A}}\mathcal{J}(\sigma),\quad \mathcal{J}(\sigma):= \sum_{n=1}^{N}\mathfrak{D}\left(\Lambda(\sigma)j_{n}, \phi_{n}\right)  
    \end{gathered}
\end{equation}
$(j_{1}, \phi_{1}),\cdots, (j_{N}, \phi_{N})$ are measurements of NtD map. The misfit function $\mathfrak{D}(u,\phi)$ measures the difference between $u$ and $\phi$. The admissible set of $\sigma$ is
\begin{equation}\label{admissible}
    \mathcal{A}=\left\{\sigma \in L_{\infty}(\Omega): c_{0} \leq \sigma \leq c_{1} \text { on } \Omega,\left.\sigma\right|_{\partial \Omega}=\left.\sigma_{0}\right|_{\partial \Omega}\right\}.
\end{equation}

In the existing method, $\mathfrak{D}$ is chosen to be the $L^{2}$ norm on $\partial \Omega$ and the regularization term $\mathcal{R}(\sigma)$ is added to $\mathcal{J}(\sigma)$ to get a new objective functional 
\begin{equation}\label{objective}
    \Psi(\sigma):=\mathcal{J}(\sigma)+\beta \mathcal{R}(\sigma),
\end{equation}
where $\beta$ is the regularization parameter. The resulting optimization problem is usually solved using iterative gradient-based optimization methods. However, due to the ill-posed nature of EIT\cite{uhlmann2009electrical}, the reconstruction is easily disturbed by the noise in the boundary measurements. A proper value of $\beta$ is required to stabilize the reconstruction process. However, when the noise reaches certain threshold level, it is hard to choose an appropriate $\beta$ to balance the smoothness and the accuracy of the reconstruction. 

Considering the favorable properties of the $W_{2}$ distance \cite{engquist2020quadratic}, here we apply $W_{2}$ for solving the optimization problem (\ref{optimization}). In fact, under the quadratic Wasserstein distance, the difference between the initial data and the disturbed data is quite small because the local cancellation of the mass makes the optimal map close to an identity map, leading to the robustness against high-frequency noise. For example \cite{villani2021topics}, $W_{2}(f_{n},f)= O(\frac{1}{n})$ and $\|f_{n}-f\|_{L^{2}}=O(1)$, where $f_{n}=1+ \sin(2\pi nx)$ and $f=1$ defined on $[0,1]$. Also, while $L^{2}$ favors the displacement along the amplitude axis, $W_{2}$ takes both spatial and amplitude changes into account. It means that $W_{2}$ is more sensitive to the shape variation of data in the space than $L^{2}$. According to Remark \ref{Monge}, the regularity theorem implies that the gradient of $W_{2}$ is two-order smoother than the input data $f$, resulting in a smoothing effect on the inversion.

To perform the optimization based on $W_{2}$, it is necessary to calculate $\mathcal{J}(\sigma)$ and $\mathcal{J}^{\prime}(\sigma)$, which involves the computation of the value and the gradient of $W_{2}$.

In the EIT problem, the misfit function $\mathfrak{D}$ is chosen to be 
\begin{equation}\label{W2}
    \mathfrak{D}(u,\phi) = W_{2}^{2}(\mathcal{L}(u),\mathcal{L}(\phi)),
\end{equation}
where $\mathcal{L}$ is the normalization operator transforming electrical potentials to non-negative density functions with unit mass, i.e. $\mathcal{L}(\phi)\in D$ in (\ref{pset}). For the EIT problem, since we measure the electrical potential $\phi\in \tilde{H}^{\frac{1}{2}}(\partial \Omega)$, the mass of $\phi$ on $\partial \Omega$ is constant zero, which satisfies mass conservation automatically. Thus we only need to rescale $\phi$ to make it positive. As mentioned in \cite{yang2018application}, a simple way is to add some positive constant $a$:
\begin{equation}\label{normalization}
    \mathcal{L}(\phi) = \frac{\phi +a}{\int_{\partial \Omega}(\phi +a )\mathrm{d}s}= \frac{1}{a}\phi +1,\quad \phi\in \tilde{H}^{\frac{1}{2}}(\partial \Omega)
\end{equation} 
In \cite{yang2018application}, the convexity property and the metric structure of $W_{2}$ are lost because of the normalization of the mass. However, we know from (\ref{normalization}) that both properties are maintained for $W_{2}(\mathcal{L}(u), \mathcal{L}(\phi))$ defined on $\tilde{H}^{\frac{1}{2}}(\partial \Omega)$, due to the fact that the mass is conserved for functions in $\tilde{H}^{\frac{1}{2}}(\partial \Omega)$. The computation of $W_{2}(\mathcal{L}(u),\mathcal{L}(\phi))$ follows from the method in the previous sections.

Another important issue is to derive the Fréchet gradient of the objective function $\mathcal{J}(\sigma)$. To simplify the notation, we discuss the case when $N=1$ in (\ref{optimization}). The first-order perturbation gives:
$$
\begin{gathered}
    \delta \mathcal{J} = \left\langle \frac{\delta \mathcal{J}}{\delta u}, \delta u\right\rangle_{\partial \Omega} = \left\langle \frac{\delta \mathfrak{D}(u,\phi)}{\delta u}, \delta u\right\rangle_{\partial \Omega}=
    \left\langle \frac{\delta \mathfrak{D}(u,\phi)}{\delta u}, \frac{\delta \Lambda(\sigma)j}{\delta \sigma}\delta \sigma\right\rangle_{\partial \Omega} \\ 
    =\left\langle \left(\frac{\delta \Lambda(\sigma)j}{\delta \sigma}\right)^{*}\frac{\delta \mathfrak{D}(u,\phi)}{\delta u}, \delta \sigma\right\rangle_{ \Omega},
\end{gathered}
$$
where $u = \Lambda(\sigma)j$ is the state variable. Thus the gradient of the functional is given by 
\begin{equation}\label{Jgradient}
    \mathcal{J}^{\prime}(\sigma) = \left(\frac{\delta \Lambda(\sigma)j}{\delta \sigma}\right)^{*}\frac{\delta\mathfrak{D}(u,\phi)}{\delta u},
\end{equation}
where the gradient 
\begin{equation}\label{D_gradient}
    \frac{\delta \mathfrak{D}(u,\phi)}{\delta u} = \frac{1}{a}\frac{\delta W_{2}^{2}(\mathcal{L}(u),\mathcal{L}(\phi))}{\delta \mathcal{L}(u)}.
\end{equation}
From Remark \ref{gradient_circle}, the gradient may be computed following the method in Section 4. For the adjoint operator of $\frac{\delta \Lambda(\sigma)j}{\delta\sigma}$, it is given in \cite{borcea2002electrical}:
\begin{equation}\label{adjoint}
    \begin{aligned}
        \left(\frac{\delta \Lambda(\sigma)j}{\delta \sigma}\right)^{*}: H^{-\frac{1}{2}}(\partial \Omega)  &\longrightarrow  L_{1}(\Omega) \\
    h\;\; & \longmapsto  -\nabla \tilde{u} \cdot \nabla u
    \end{aligned}
\end{equation}
where $u= \Lambda(\sigma)j$ and $\tilde{u}=\Lambda(\sigma)h$. Note that the Fréchet gradient of $W_{2}$ in (\ref{D_gradient}) involves the choice of a hyper-parameter $c$, as stated in (\ref{fre_exp}). Here the constant $c$ is chosen to make 
$$\int_{\partial \Omega} \frac{\delta \mathfrak{D}(u,\phi)}{\delta u}\mathrm{d}s =\int_{\partial \Omega}\frac{\delta W_{2}^{2}(\mathcal{L}(u),\mathcal{L}(\phi))}{\delta u} \mathrm{d}s=0, $$
hence $\frac{\delta \mathfrak{D}(u,\phi)}{\delta u} \in \tilde{H}^{-\frac{1}{2}}(\partial \Omega)$, which makes it possible to apply the operator (\ref{adjoint}) to this gradient.

With the gradient information, we propose the Barzilai-Borwein gradient algorithm with a non-monotone line search strategy to minimize (\ref{objective}). The detailed implementation of the method is summarized in Algorithm \ref{alg2}.
To make the optimization procedure more stable, following \cite{jin2012reconstruction}, \cite{knowles1998variational}, we apply a preconditioned gradient of $\mathcal{J}(\sigma)$, namely the sobolev gradient $\mathcal{J}_{s}^{\prime} := (I-\Delta)\mathcal{J}^{\prime}$, which solves the equation \cite{neuberger2009sobolev}:
\begin{equation}\label{Jsgradient}
\begin{aligned}
-\Delta \mathcal{J}_{s}^{\prime}(\sigma)+\mathcal{J}_{s}^{\prime}(\sigma)& =\mathcal{J}^{\prime}(\sigma)\; \text { in } \Omega \\
\mathcal{J}_{s}^{\prime}(\sigma) &=0 \;\qquad \text { on } \partial \Omega.
\end{aligned}
\end{equation}
Since $\mathcal{J}_{s}^{\prime}$ is zero on $\partial \Omega$, it naturally satisfies the boundary condition in (\ref{admissible}).

\begin{algorithm}  
    \caption{The Barzilai-Borwein gradient method}   
    \label{alg2}  
    \begin{algorithmic}  
    \REQUIRE initial $\sigma_{0}$, integer $M\geq 0$, $\tau\in(0,1)$, $0<\rho_{1}<\rho_{2}<1$, $s_{min}$ and $s_{max}$; stop criterion $s_{stop}$ and $I_{max}$.  
    \STATE Set $k:=0$, $s_{0} := s_{max}$, $\sigma_{+} := \sigma_{0}$, and compute $\mathcal{J}_{s}^{\prime}(\sigma_{0})$.
    \WHILE{$s_{k}> s_{stop}$ and $k < I_{max}$}
    \WHILE{\[\Phi(\sigma_{+})\geq \max_{0\leq j\leq M-1}\Phi(\sigma_{k-j}) - \frac{\tau}{2s_{k}} \|\sigma_{+}- \sigma_{k}\|_{H^{1}(\Omega)}^{2},\]}
    \STATE Select $\rho\in[\rho_{1}, \rho_{2}]$ and update $s_{k}:=\rho *s_{k}$.
    \STATE Let $\gamma_{k} = \sigma_{k} - s_{k}\mathcal{J}_{s}^{\prime}(\sigma_{k})$. Solving the proxy problem
    $$\sigma_{+} := \mathop{\mathrm{argmin}}\limits_{\sigma\in\mathcal{A}}\, \frac{1}{2s_{k}}\left\|\sigma - \gamma_{k}\right\|_{H^{1}(\Omega)}^{2} + \beta \mathcal{R}(\sigma).$$
    \ENDWHILE
    \STATE Let $\sigma_{k+1}:=\sigma_{+}$ and compute $\mathcal{J}_{s}^{\prime}(\sigma_{k+1})$. Then calculate 
    $$
    \begin{gathered}
        x_{k}:= \langle \sigma_{k+1}-\sigma_{k},\sigma_{k+1}-\sigma_{k}\rangle_{H_{1}(\Omega)}\\
        y_{k}:= \left\langle \sigma_{k+1}-\sigma_{k} ,\mathcal{J}_{s}^{\prime}(\sigma_{k+1}) - \mathcal{J}_{s}^{\prime}(\sigma_{k})\right\rangle_{H_{1}(\Omega)}.
    \end{gathered}
    $$
    \IF{$y_{k}\leq 0$}
    \STATE $s_{k+1}:=s_{max}$
    \ELSE
    \STATE $s_{k+1}:= \min\{\,s_{max},\max\{\,s_{min}, \frac{x_{k}}{y_{k}}\,\}\,\}$. 
    \ENDIF
    \STATE $k:=k+1$
    \ENDWHILE
    \STATE The output $\sigma_{k}$ is the final reconstructed result.
    \end{algorithmic}  
\end{algorithm}
\section{Numerical Experiments}
In the following experiments, the forward problems and adjoint problems are solved numerically using the finite element method in FEniCS \cite{logg2012automated}. The current data set is chosen to be
\begin{equation}\label{set}
\left\{\, \sin(n\theta),\, \cos(n\theta)\,|\,n=1,\cdots,N\,\right\}
\end{equation}
A Gaussian noise is added to the measured data, i.e., the electrical potentials take
$$\phi_{n} := \phi_{n} + \xi *\varepsilon\max_{k}\|\phi_{k}\|_{L^{\infty}(\partial \Omega)},\quad \xi \sim \mathcal{N}\left(0, 1\right)$$
where $\varepsilon$ refers to the relative noise level.

\begin{figure*}[h]
    \centering
    \includegraphics[width=0.4\textwidth]{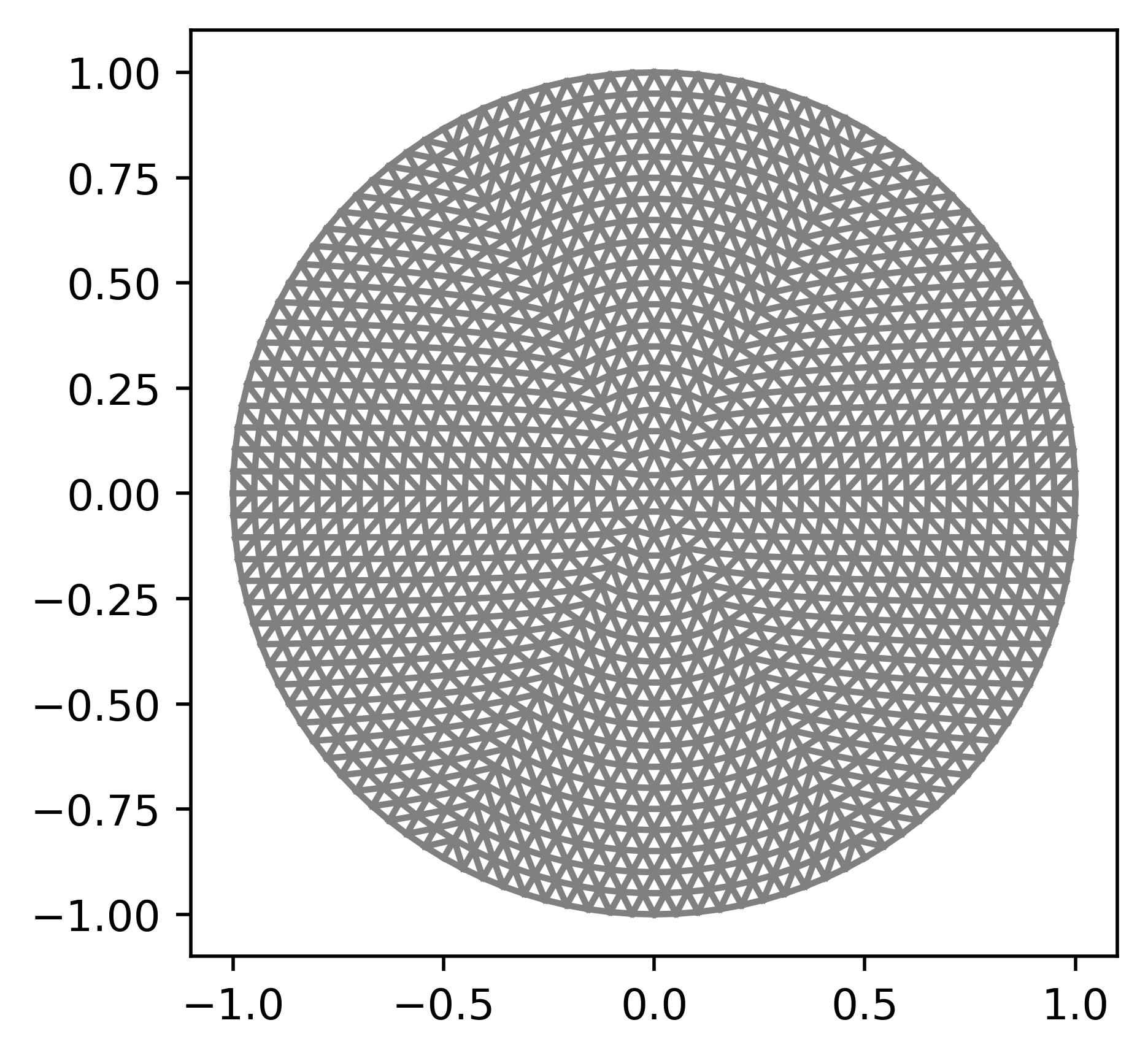}
    \caption{Mesh for solving the inverse problem.}\label{mesh}
\end{figure*}

The inverse problem is discretized using piecewise linear finite elements with 2400 triangle elements (Figure \ref{mesh}), and the exact data is calculated on a refined mesh.
In the following, $N$ in (\ref{set}) is set to be five, which means ten NtD data are generated to solve the inverse problem. In algorithm, the minimal step size $s_{min}:= 1$ and the maximum step size $s_{max}:=1000$. The parameters $M:=5$, $\tau := 1\times 10^{-5}$, $\rho_{1} :=0.4$, $\rho_{2} := 0.6$ are also fixed throughout all the experiments. The stopping criteria are the lower bound on the step size $s_{stop}:=1\times 10^{-3}$ and the maximum number of iterations $I_{max}$. It should be noted that the model will fail to converge when high-level noise is applied, and thus $s_{stop}$ is invalid. In order to give feasible reconstructions, $I_{max}$ is set to be different in each experiment. 

In the meantime, we also give the reconstructions based on the conventional $L^{2}$ norm for comparison. Different examples are presented to illustrate the features of $W_{2}$ inversion.  Without specification, the relative noise level $\varepsilon = 3\%$ and the initial conductivity $\sigma_{0}$ is chosen to be constant $1$. 

\begin{example}\label{6_1}
    We start with a simple example. The expression of the true conductivity field $\sigma$ (Figure \ref{1_true}) is given by
\begin{equation}
    \sigma(x, y)=\left\{\begin{array}{ll}2, & (x+0.3)^{2}+(y-0.3)^{2} \leq (0.35)^{2} \\ 1, & \text { else }\end{array}\right.
\end{equation}
\begin{figure*}[h]
    \centering
    \subfloat[True $\sigma$]{\includegraphics[width=0.32\textwidth]{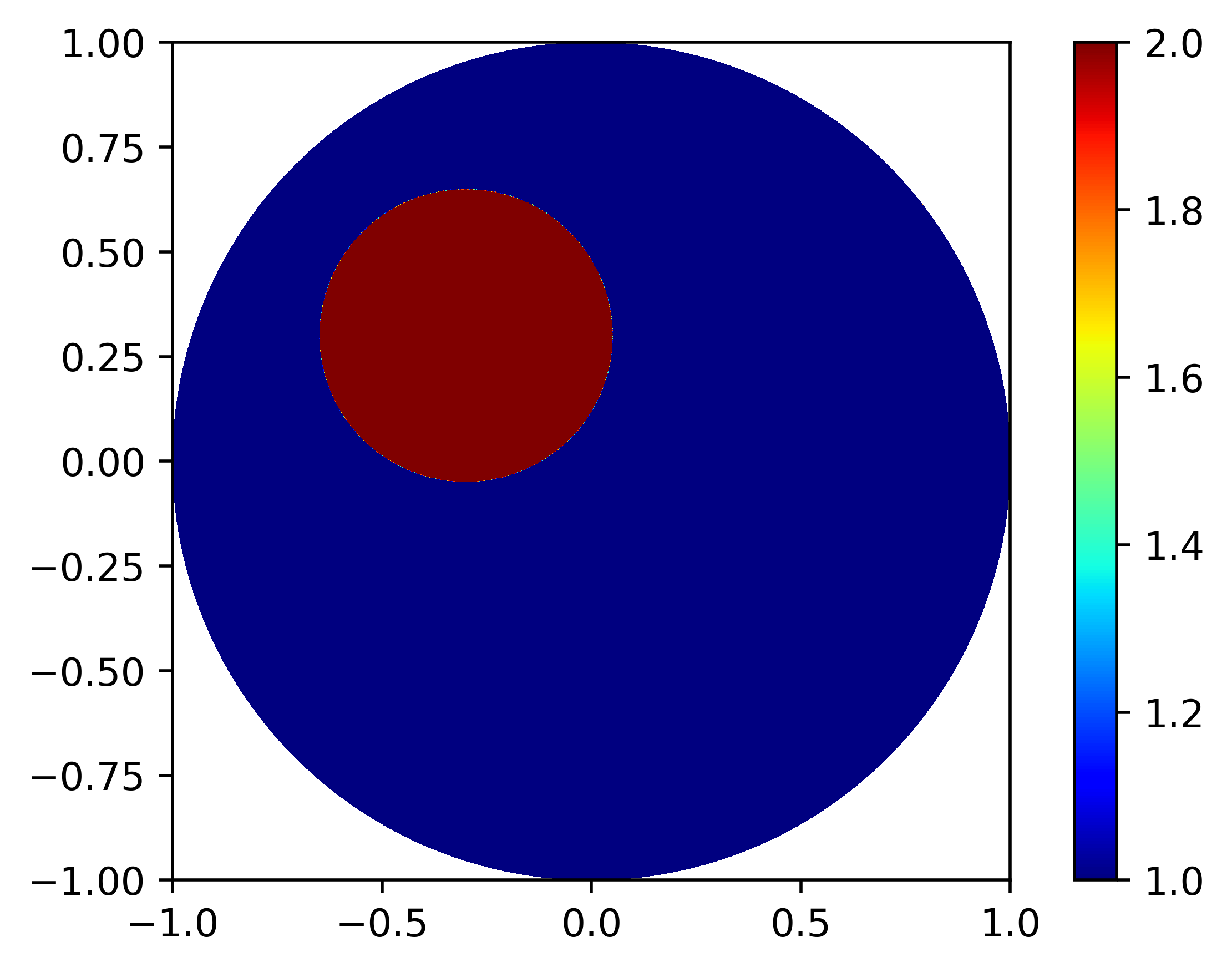}\label{1_true}}
    \;
    \subfloat[$W_{2}$]{\includegraphics[width=0.32\textwidth]{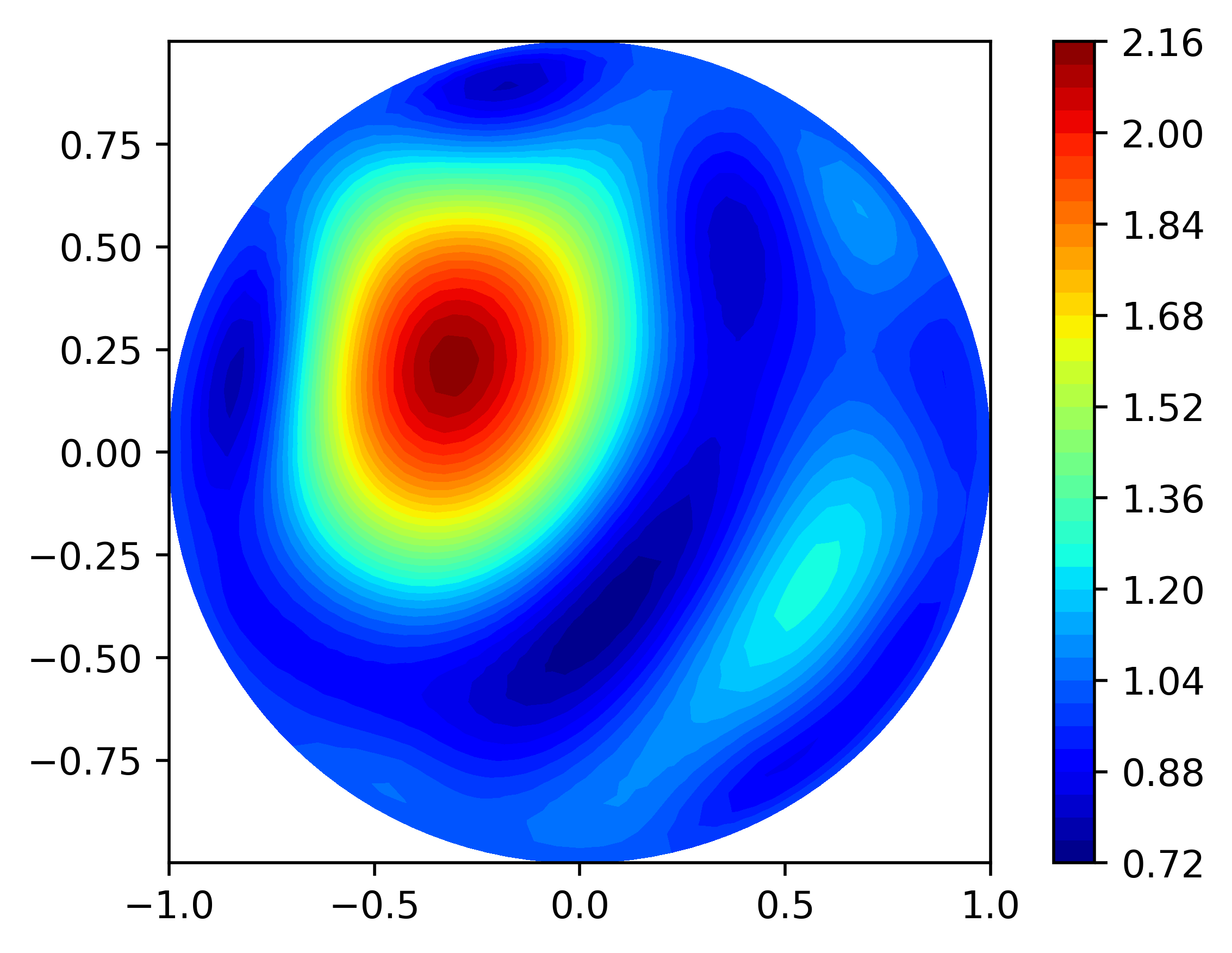}\label{1_W}}
    \;
    \subfloat[$L^{2}$]{\includegraphics[width=0.32\textwidth]{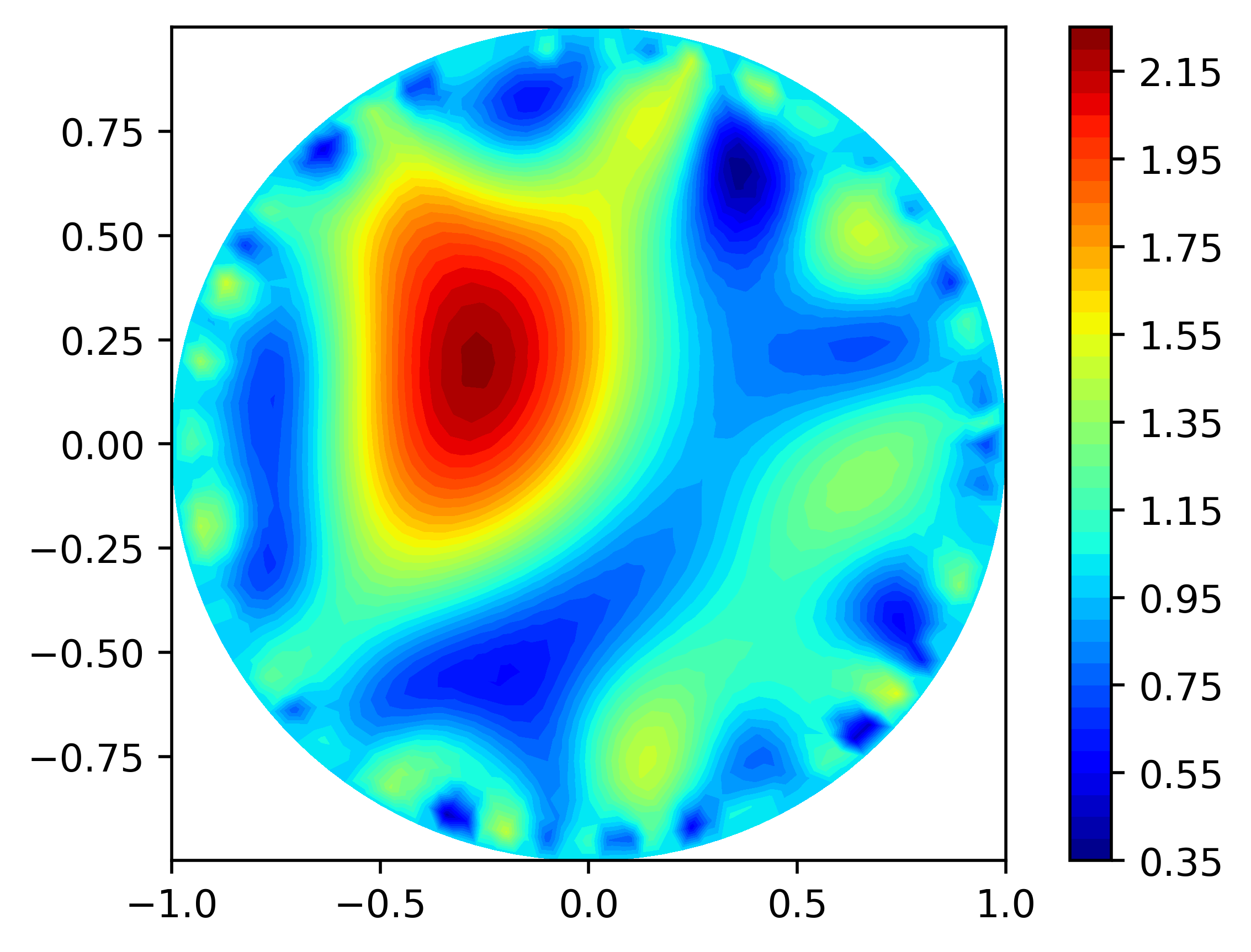}\label{1_L}}
    \caption{Results of Example \ref{6_1} with $3\%$ in the data. }\label{fig1}
  \end{figure*}
The regularization parameter $\beta$ is set to be zero for both methods.

Figure \ref{1_W} and \ref{1_L} show the reconstruction results by the $W_{2}$ distance and $L^{2}$ distance after 500 iterations, respectively. As shown in \cite{isaacson1986distinguishability}, the NtD measurements are insensitive to the magnitude of the inclusion conductivity. In most situations, it is hard to evaluate the magnitude since instabilities would take over after a modest number of iterations, especially with noise present \cite{knowles1998variational}. This phenomenon will also be demonstrated in Example \ref{6_2} and \ref{6_3}. To reveal all the information in the data, we set $I_{max}=500$. As shown in Figure \ref{fig1}, with slightly perturbed by the noise, the inclusion's magnitude, shape, and position are basically retrieved by both methods. However, while the image of $L^{2}$ inversion is heavily polluted by noise, $W_{2}$ inversion shows a strong resilience in the presence of high-level noise and a huge number of iterations. 

\end{example}

\begin{example}\label{6_2}
The inclusion conductivity is set to be an ellipse (Figure \ref{2_true}):
$$\sigma(x, y)=\left\{\begin{array}{ll}2, & \frac{x^{2}}{0.04}+\frac{(y-0.5)^{2}}{0.16} \leq 1 \\ 1, & \text { else }\end{array}\right.$$

\begin{figure*}[h]
    \centering
    \includegraphics[width=0.4\textwidth]{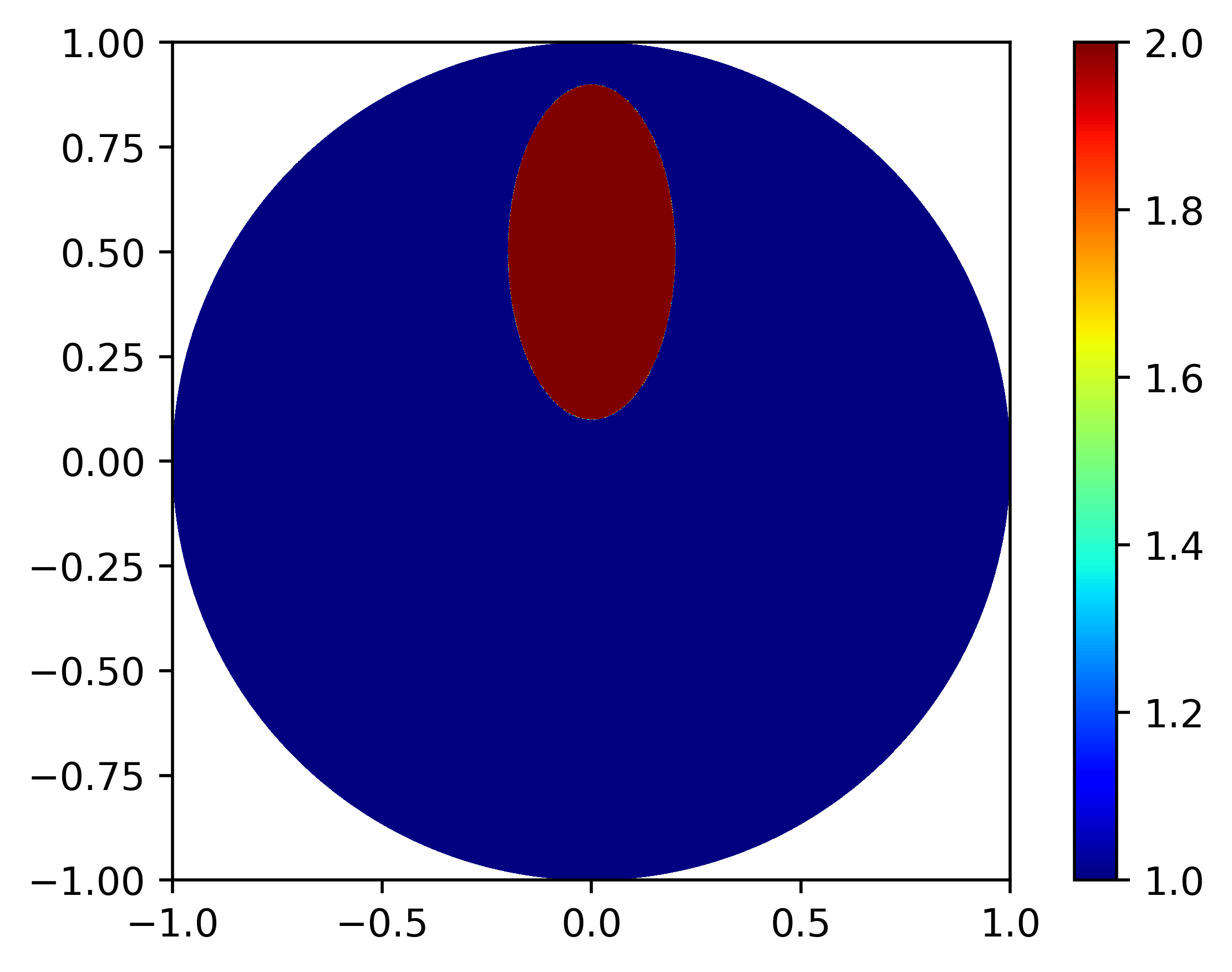}
    \caption{True conductivity field of Example \ref{6_2}}\label{2_true}
\end{figure*}

We choose total variation regularization scheme for $L^{2}$ inversion with $\beta=1\times 10^{-4}$. No regularization is applied for $W_{2}$ inversion. The evolutions of the reconstructions are plotted in Figure \ref{fig2} at iteration 25, 50, and 100.
\begin{figure}[h] 
    \centering
    \subfloat[$W_{2}$: 25 iter]{\includegraphics[width=0.34\textwidth]{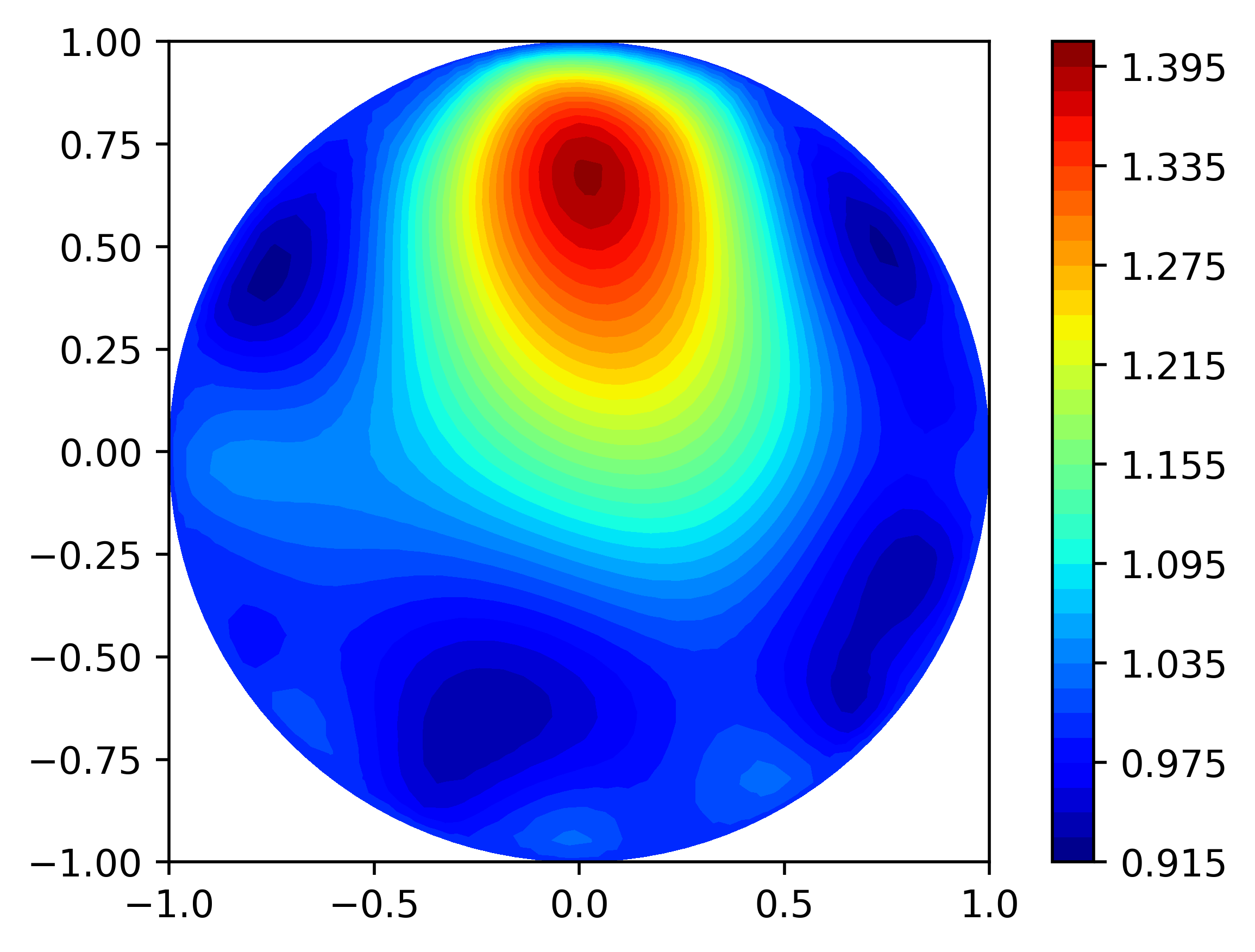}\label{W_25}}
    \subfloat[$W_{2}$: 50 iter]{\includegraphics[width=0.33\textwidth]{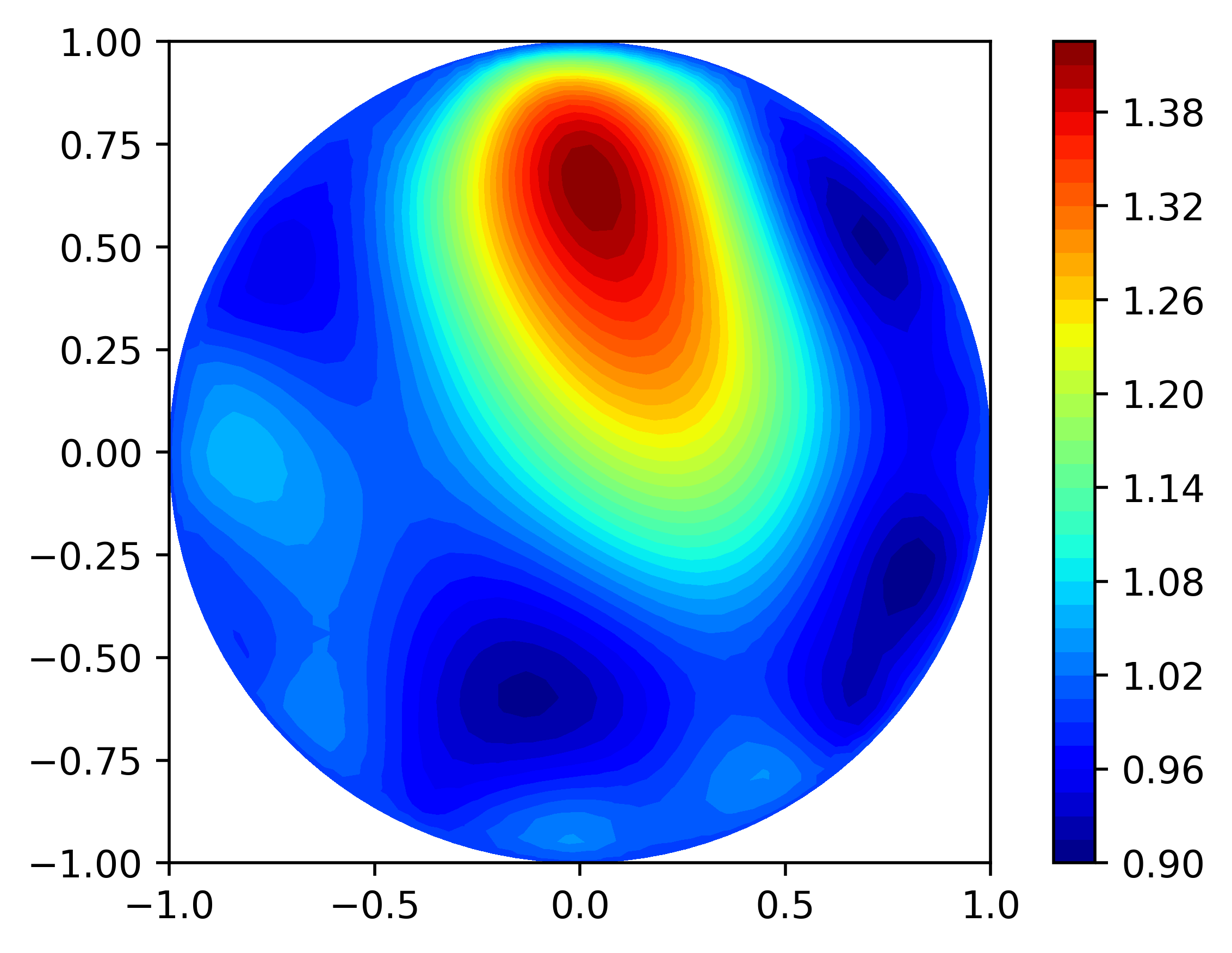}\label{W_50}}
    \subfloat[$W_{2}$: 100 iter]{\includegraphics[width=0.335\textwidth]{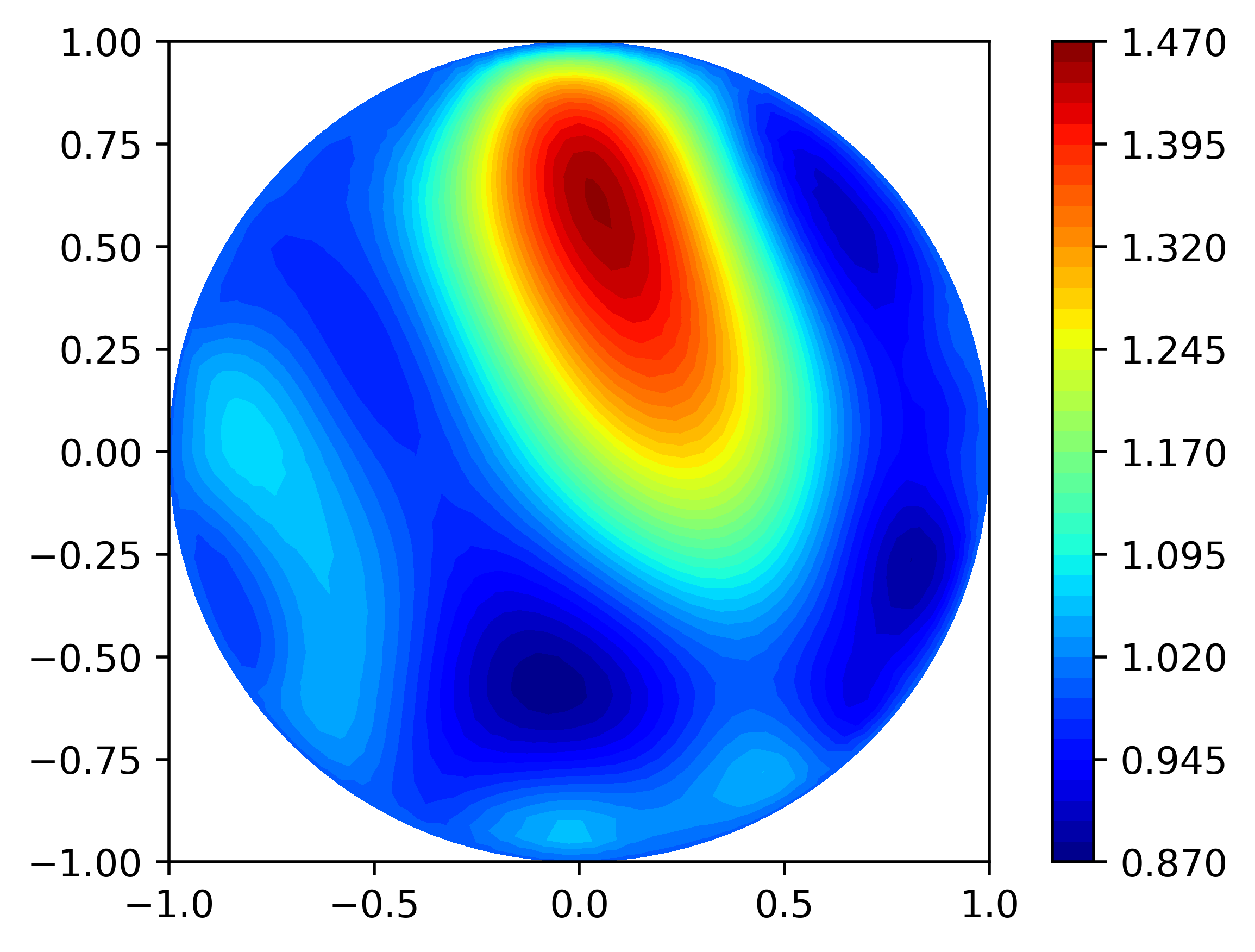}\label{W_100}}

    \subfloat[$L^{2}$: 25 iter]{\includegraphics[width=0.34\textwidth]{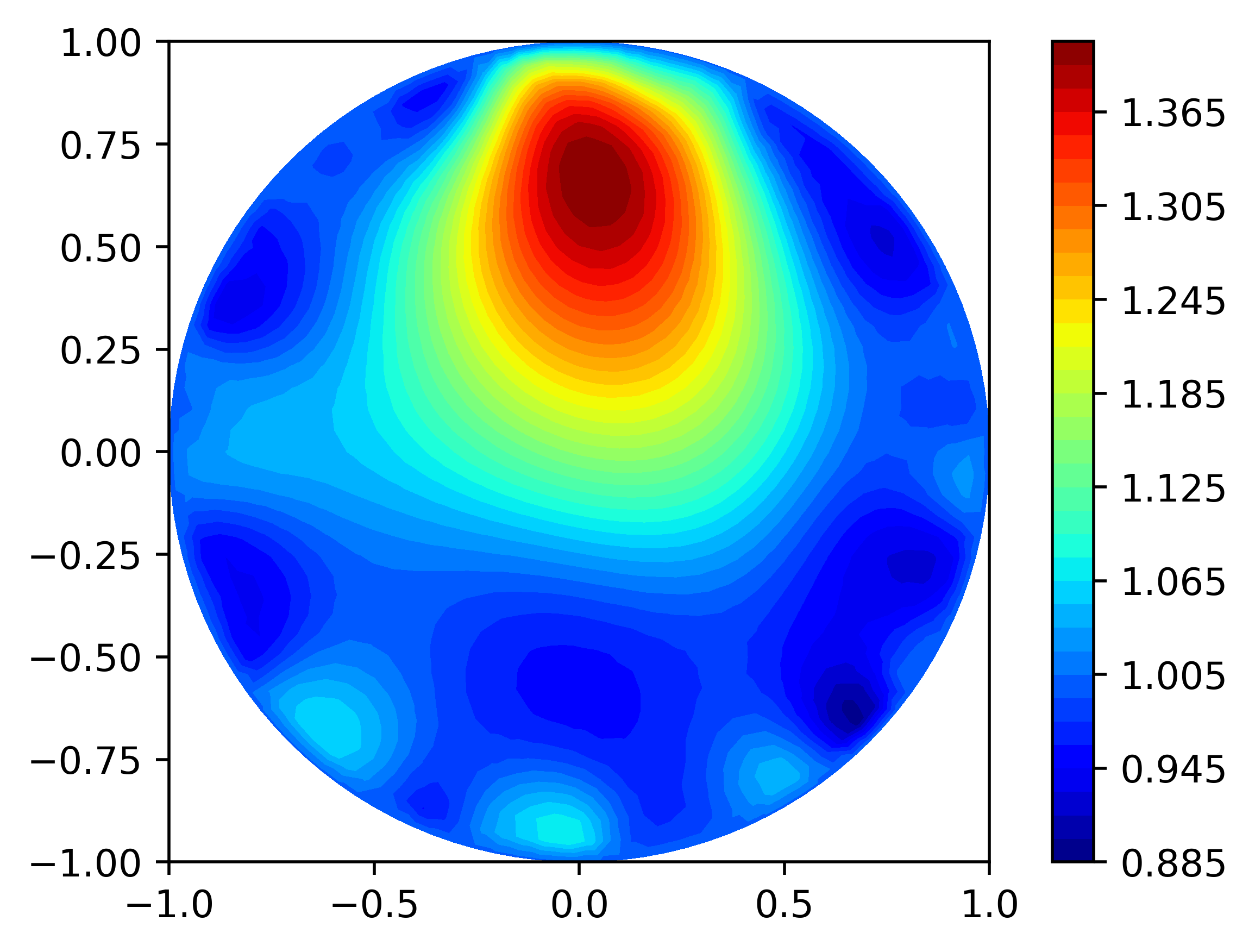}\label{L_25}}
    \subfloat[$L^{2}$: 50 iter]{\includegraphics[width=0.33\textwidth]{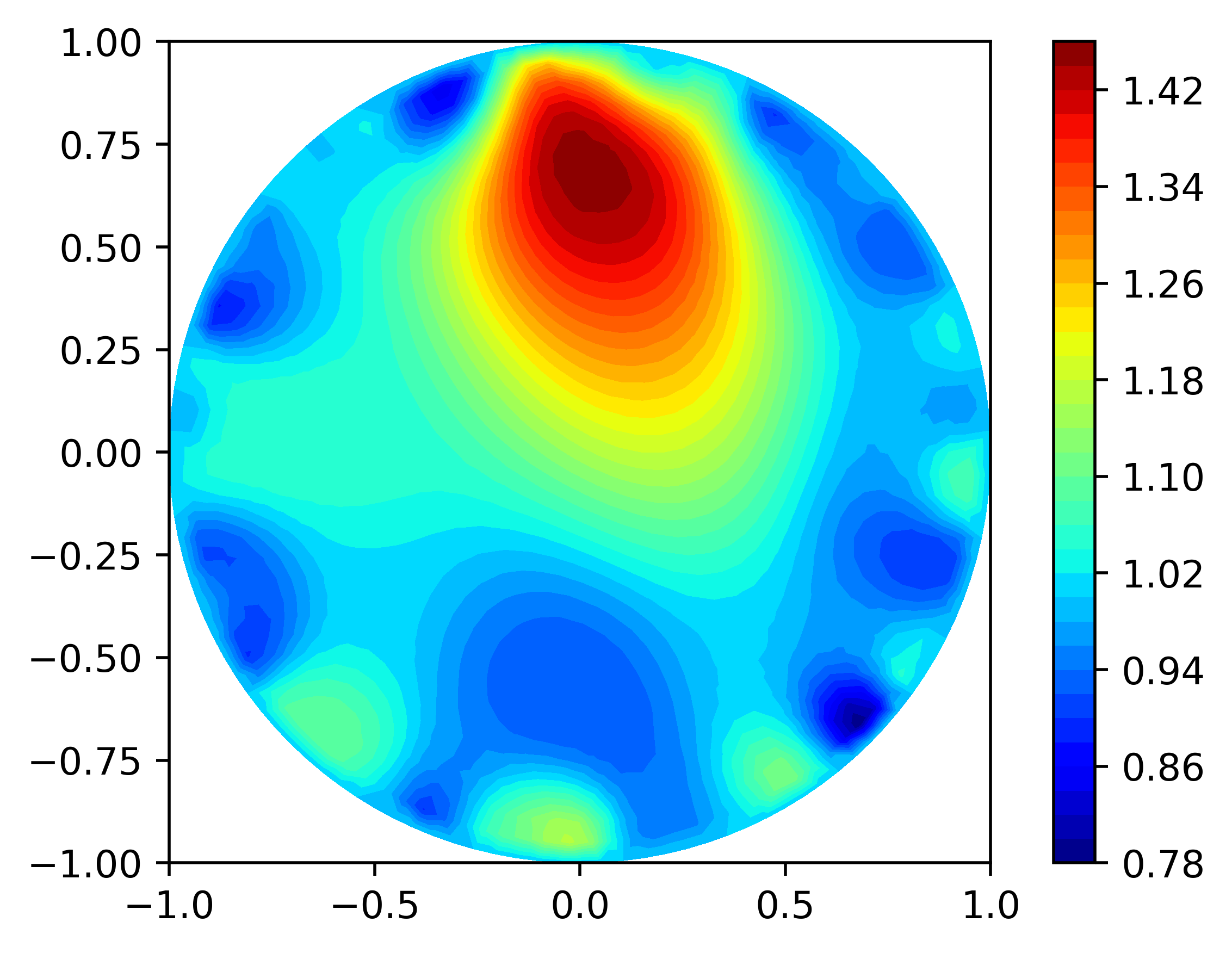}\label{L_50}}
    \subfloat[$L^{2}$: 100 iter]{\includegraphics[width=0.33\textwidth]{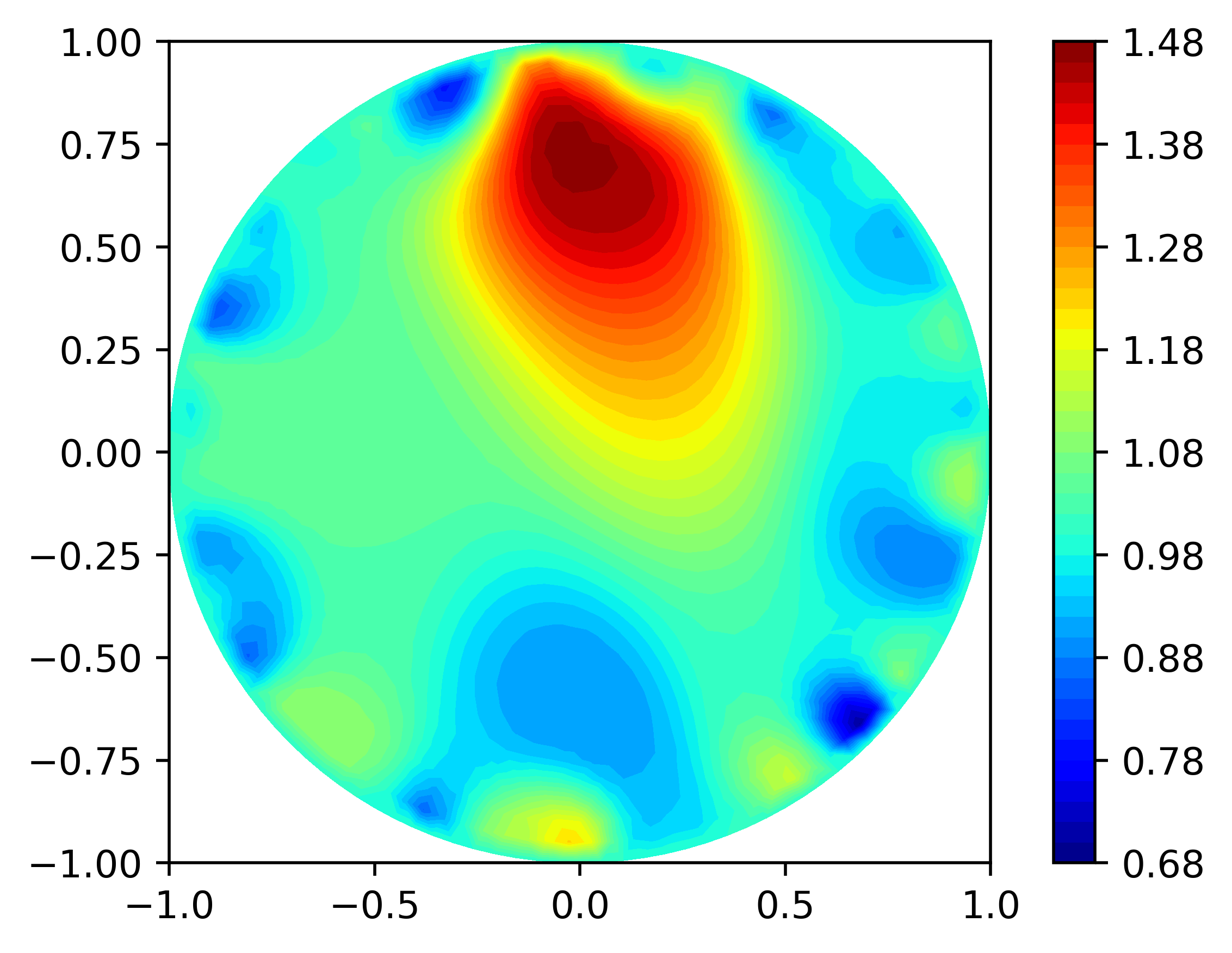}\label{L_100}}
    \caption{Reconstructions for Example \ref{6_2} at different iterations with $3\%$ noise in the data. The first row: $W_{2}$ reconstructions. The second row: $L^{2}$ reconstructions with $\beta=1\times 10^{-4}$.}\label{fig2}
\end{figure}
Compared with Example \ref{6_1}, the shape of the inclusion is much more singular in this example. As one should expect, recovery of the conductivity is easily degenerated at a high noise level while the iteration progressing. To better preserve the reconstructed shape, we stop the reconstruction process at iteration $I_{max}=100$. Hence the magnitude of the inclusion is severely underestimated. 

The reconstructions are pretty smooth and regular at lower iterations. However, information in the data can not be fully revealed, which leads to a significant error in the recovery. As the iterations increase, severe artifacts appear in the reconstructed images near the boundary and the edge around the inclusion when using the $L^{2}$ distance for inversion. What's more, the inclusion is deformed due to the high-level noise, causing the inaccuracy of the aspect ratio. In contrast, noise in the inversion result of the $W_{2}$ distance can be significantly smoothed. Despite the noise-induced rotation of the inclusion, $W_{2}$ inversion gives a fair representation of the inclusion's shape. The elliptical boundary of the inclusion is smooth and clearly visible, with an aspect ratio close to 2:1.

It should be noted that the quality of the image is difficult to be improved unless there is strong prior information about the conductivity. For example, the conductivity field is a piecewise constant with a known background. When such a priori knowledge is assumed, sparsity regularization can be employed to improve the reconstruction quality significantly.
\end{example}

\begin{example}\label{6_3}
The true conductivity field is shown in Figure \ref{3_true}, which consists of three elliptical inclusions. The magnitude of the upper one is 0.5, and the bottom two is 2. This example is to simulate a cross-section of the human chest. Multiple inclusions with higher and lower inclusions are challenging for some numerical algorithms since the envelope of the conductivity is non-convexed. 
    \begin{figure}[h]
        \centering
        \subfloat[True $\sigma$]{\includegraphics[width=0.4\textwidth]{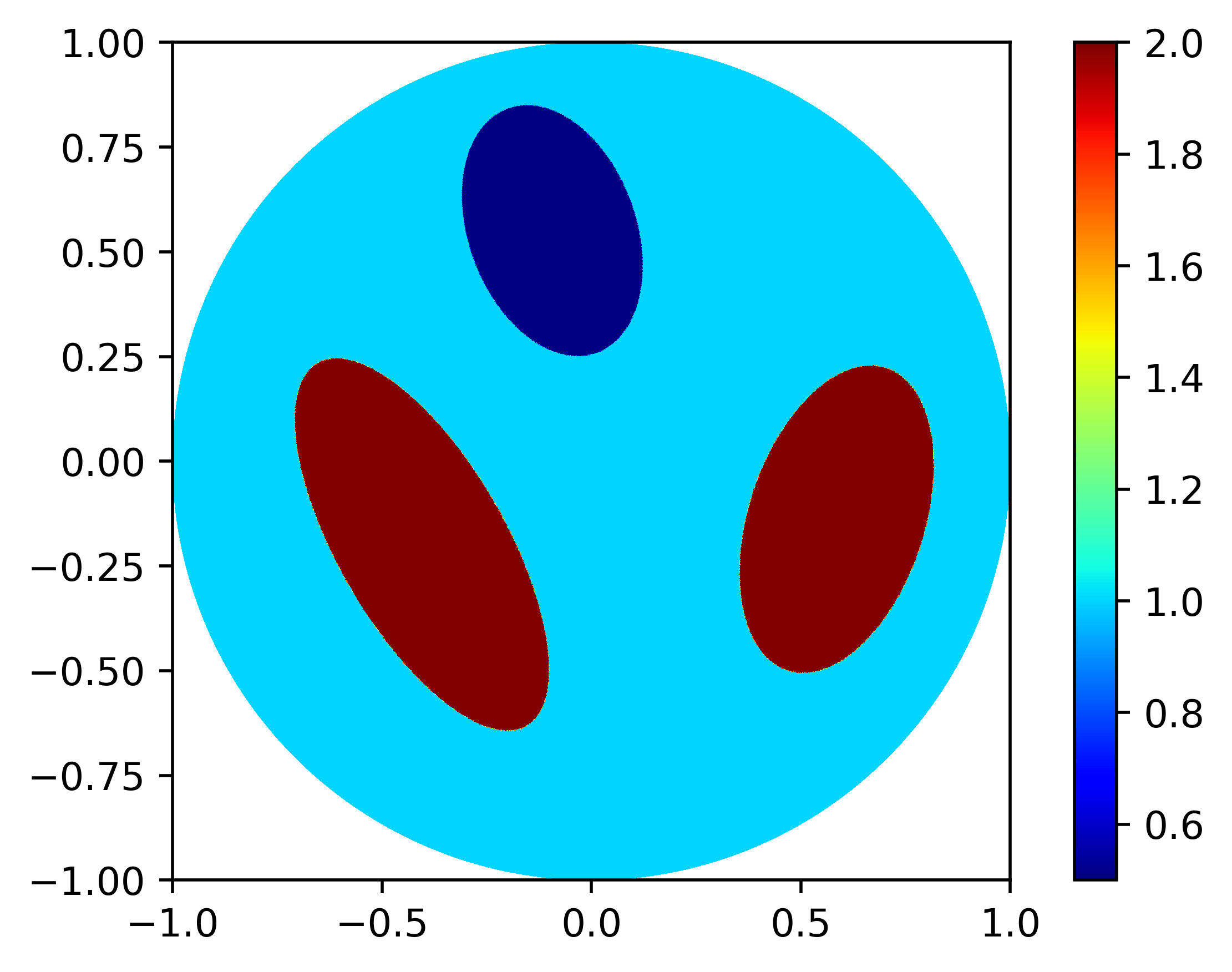}\label{3_true}}
        \subfloat[$W_{2}$]{\includegraphics[width=0.41\textwidth]{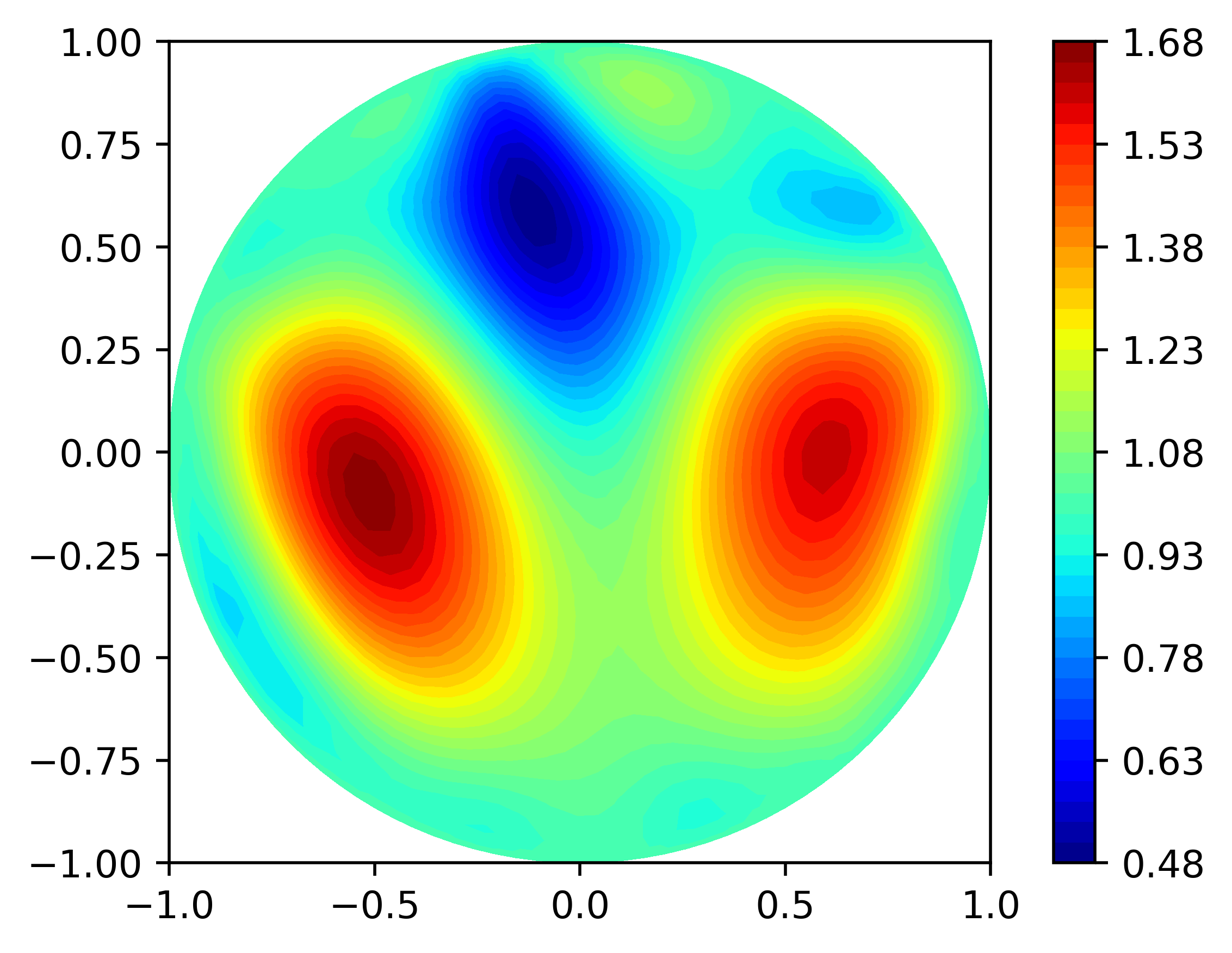}\label{3_W}}

        \subfloat[$L^{2}: \beta = 1\times 10^{-4}$]{\includegraphics[width=0.41\textwidth]{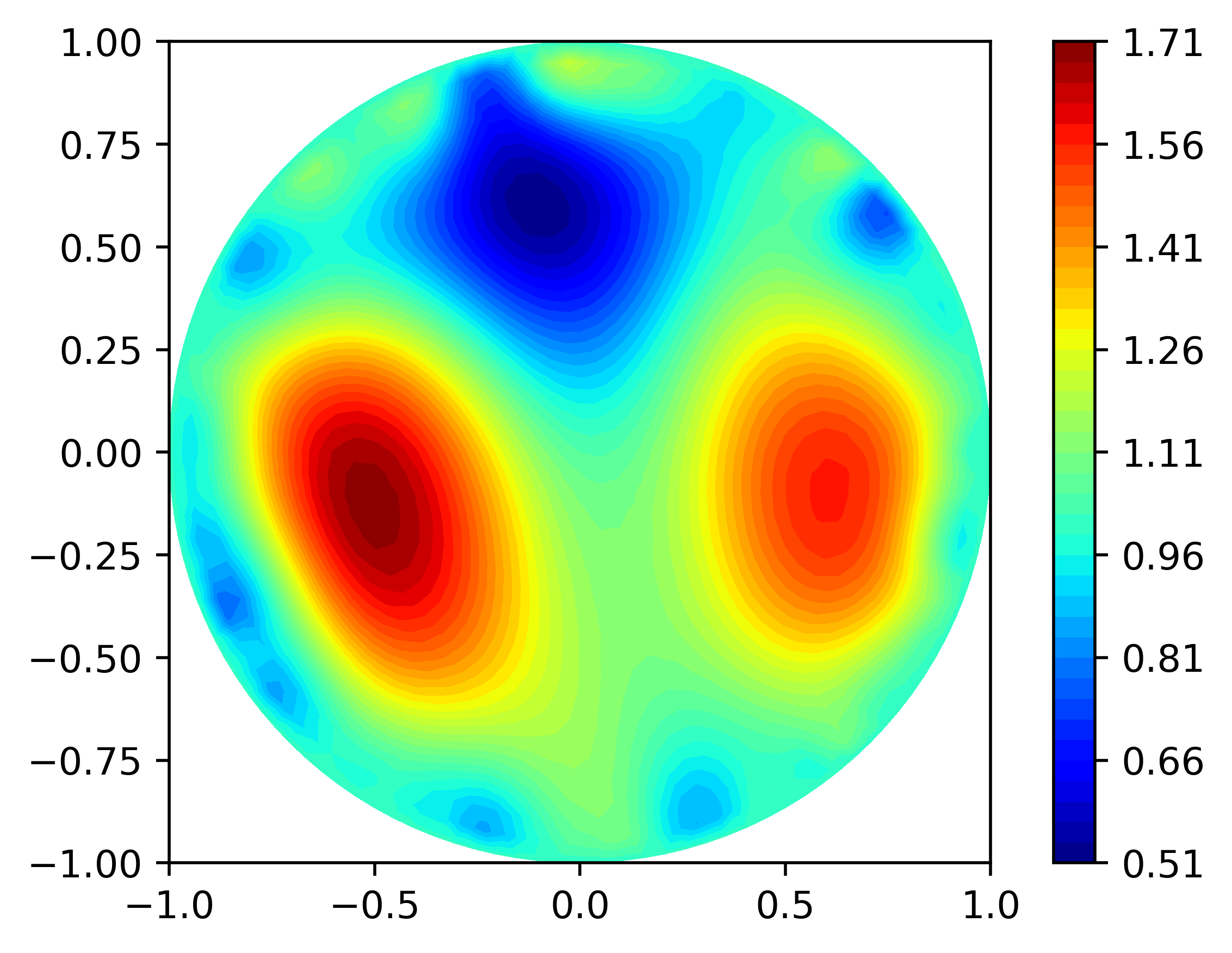}\label{3_L}}
        \subfloat[$L^{2}: \beta = 9\times 10^{-4}$]{\includegraphics[width=0.41\textwidth]{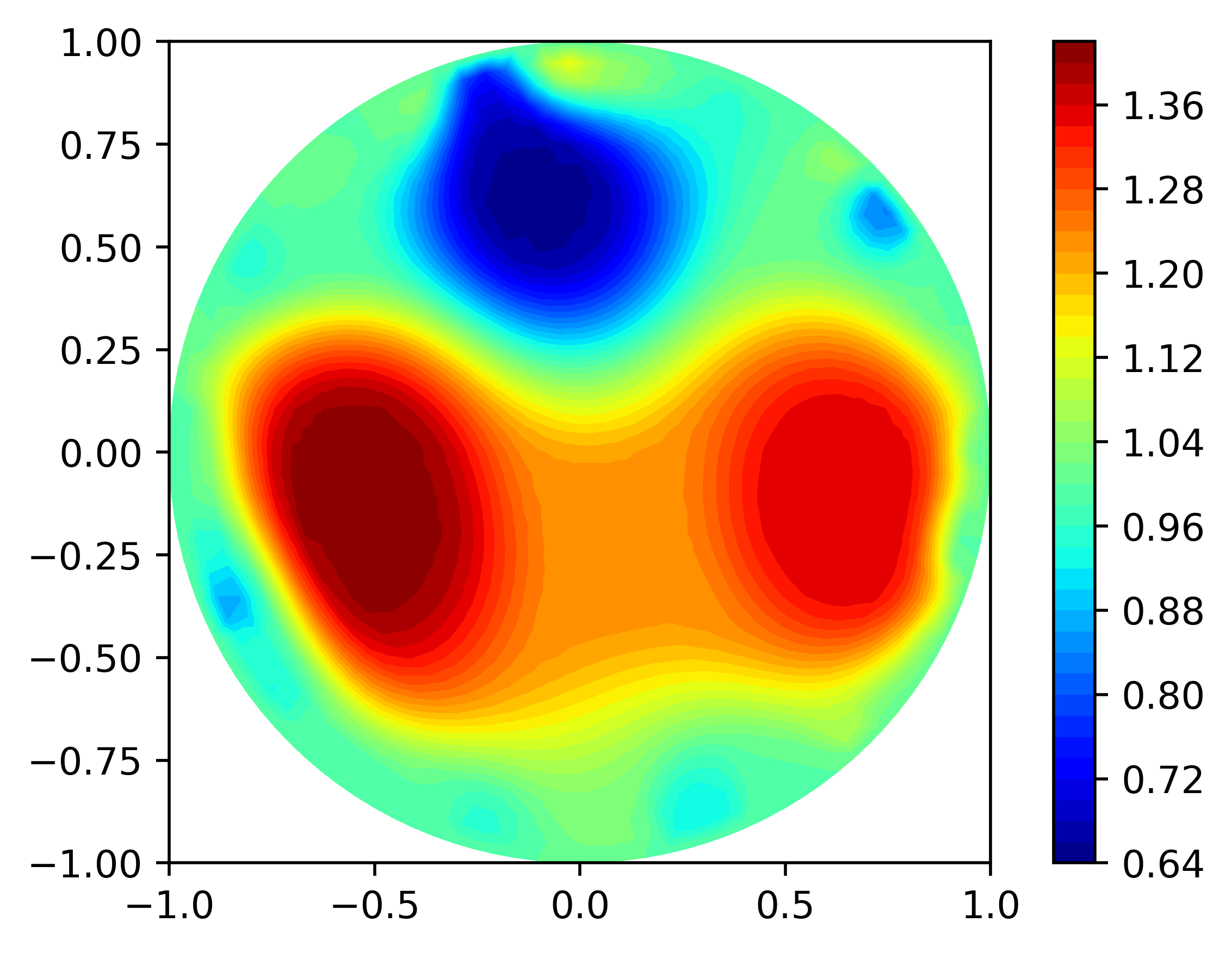}\label{3_LL}}
        \caption{Results for Example \ref{6_3} with $3\%$ noise in the data. The first row: (\ref{3_true}) true conductivity of Example \ref{6_3} and (\ref{3_W}) $W_{2}$ reconstruction after 70 iterations. The second row: (\ref{3_L}) $L^{2}$ reconstruction with $\beta=1\times 10^{-4}$ after 70 iterations and (\ref{3_LL}) $L^{2}$ reconstruction with $\beta=9\times 10^{-4}$ after 100 iterations.}\label{fig3}
    \end{figure}
Two different regularization parameters are used for the inversion with the $L^{2}$ distance. Figure \ref{fig3} shows the numerical result after appropriate iterations. In the $L^{2}$ reconstruction with $\beta=1\times 10^{-4}$, three inclusions are identified. However, the upper inclusion is severely distorted. This phenomenon is attributed to some specific properties of the upper inclusion. On the one hand, the change of magnitude is less obvious in this ellipse, and thus it is more difficult to detect its features. On the other hand, while the upper half is close to the boundary and susceptible to noise perturbation, the lower half of this ellipse is far away from the boundary, making it hard to be reconstructed.

Compared with Figure \ref{3_L}, the structure of the reconstruction in Figure \ref{3_W} is much more stable. Not only are small spurious oscillations erased, but the shape and the relative position are retrieved in a precise way.

Then we perform an experiment using the $L^{2}$ distance with a larger regularization parameter $\beta =9\times 10^{-4}$. The reconstructed image is shown in Figure \ref{3_LL}. Compared with Figure \ref{3_L}, noise near the boundary is somewhat eliminated, and the inclusions become more regular. Nonetheless, the regularization parameter is so large that the magnitude of the conductivity is significantly underestimated. The shape of the upper inclusion is still not correctly identified, and the demarcation of the inclusions is blurred. The result suggests that even large parameters fail to achieve the regularization effect of $W_{2}$ distance.
\end{example}

\begin{example}\label{6_4}
As shown in Figure \ref{fig4}, the true conductivity field is described as one single circular inclusion plus a homogeneous background. The inclusion is centered at $(0.5, \frac{3}{4}\pi)$ in the polar coordinate, with a radius $0.22$. Suppose that the size and the shape of the inclusion are already known. The inverse problem is to determine the center of the inclusion from NtD measurements at (\ref{set}) with $N=1$. $10\%$ noise is added to the observed data. No regularization is applied to the objective functions. 
\begin{figure*}[h]
    \centering
    \includegraphics[width=0.4\textwidth]{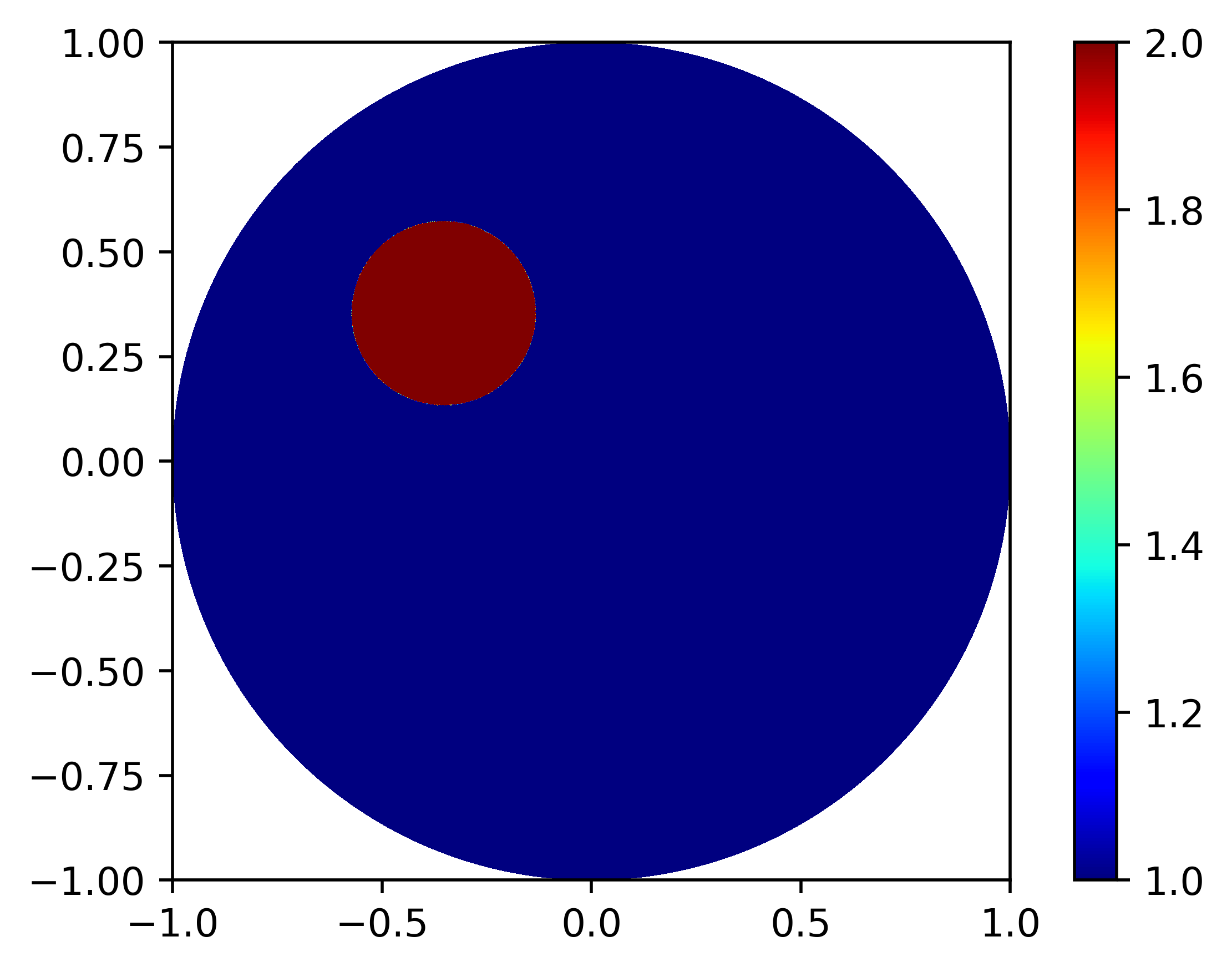}
    \caption{True conductivity field of Example \ref{6_4}.}\label{fig4}
\end{figure*}
The landscapes of the $W_{2}$ and $L^{2}$ objective functions are illustrated in Figure \ref{fig5}. Due to the effect of the noise, the $L^{2}$ landscape becomes extremely oscillatory and possesses many local minima, especially near the boundary. Nevertheless, the optimization landscape of $W_{2}$ is much smoother, which has been analyzed in \cite{engquist2020quadratic} from a local viewpoint. As we show in Figure \ref{W_comp} and \ref{L_comp}, local minima can be smoothed out by $W_{2}$ inversion. 
    \begin{figure}[h]
        \centering
        \subfloat[$W_{2}$ landscape]{\includegraphics[width=0.4\textwidth]{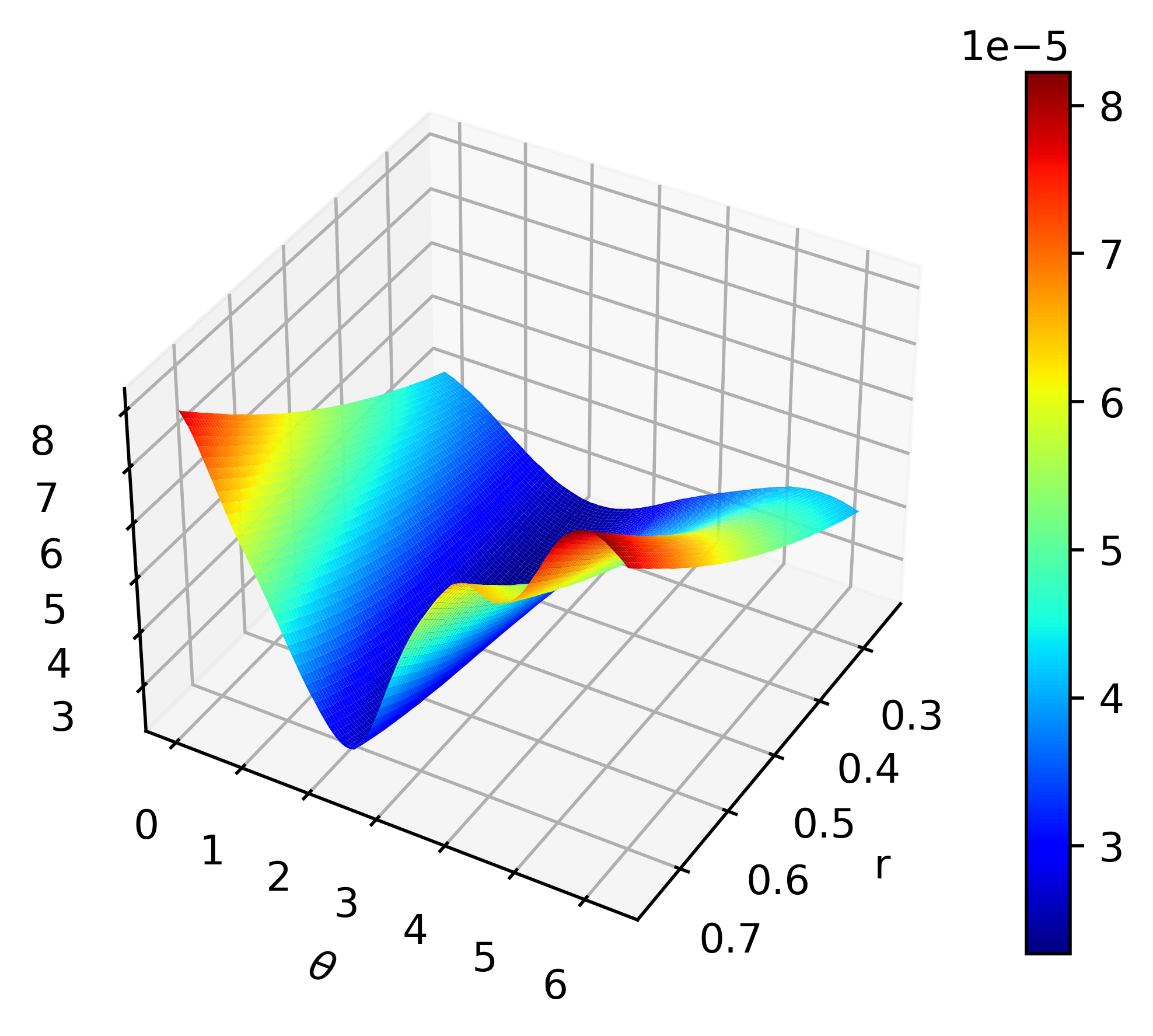}\label{W_land}}
        \quad\;
        \subfloat[$L^{2}$ landscape]{\includegraphics[width=0.45\textwidth]{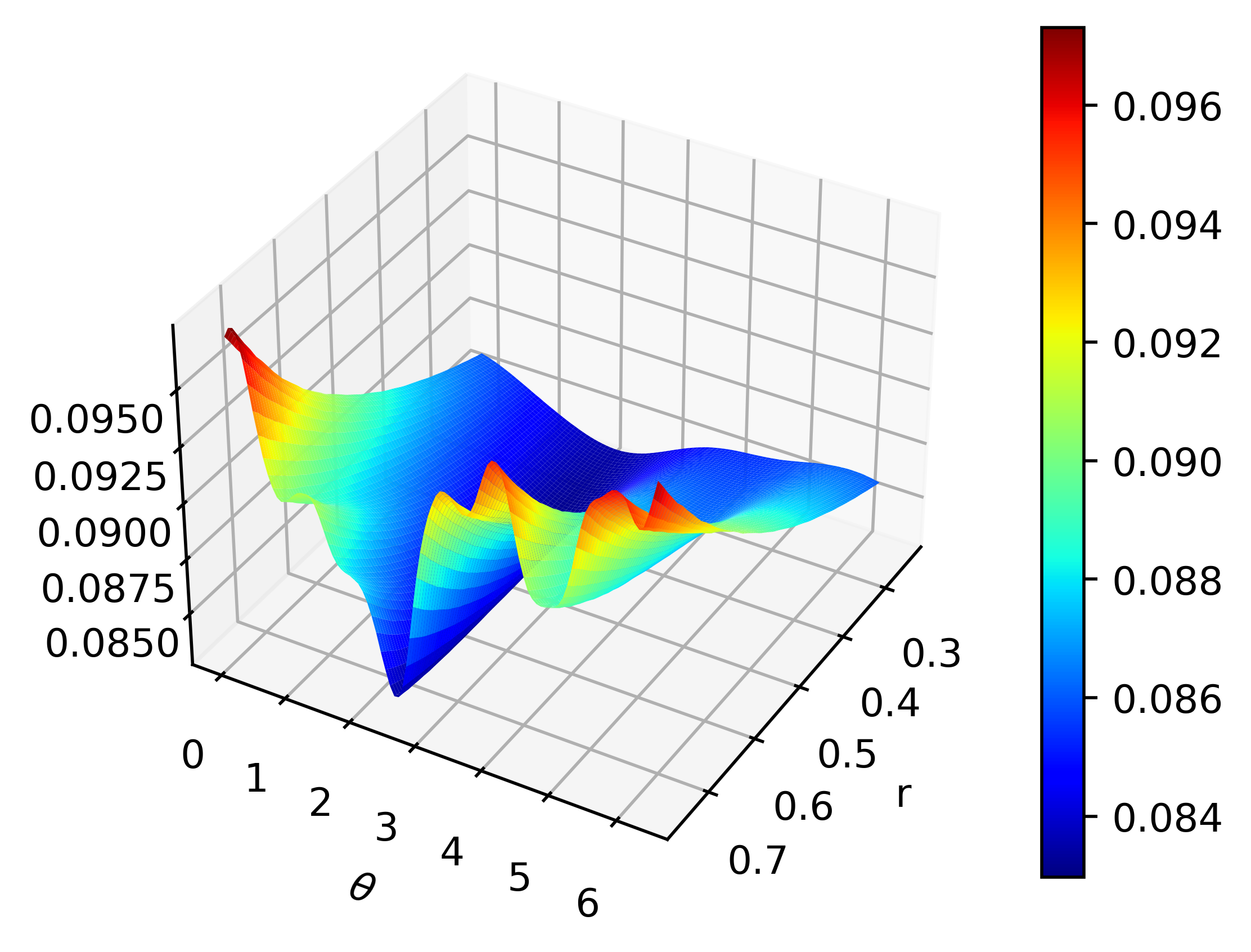}\label{L_land}}
    
        \subfloat[$W_{2}$ landscape at R = 0.5]{\includegraphics[width=0.4\textwidth]{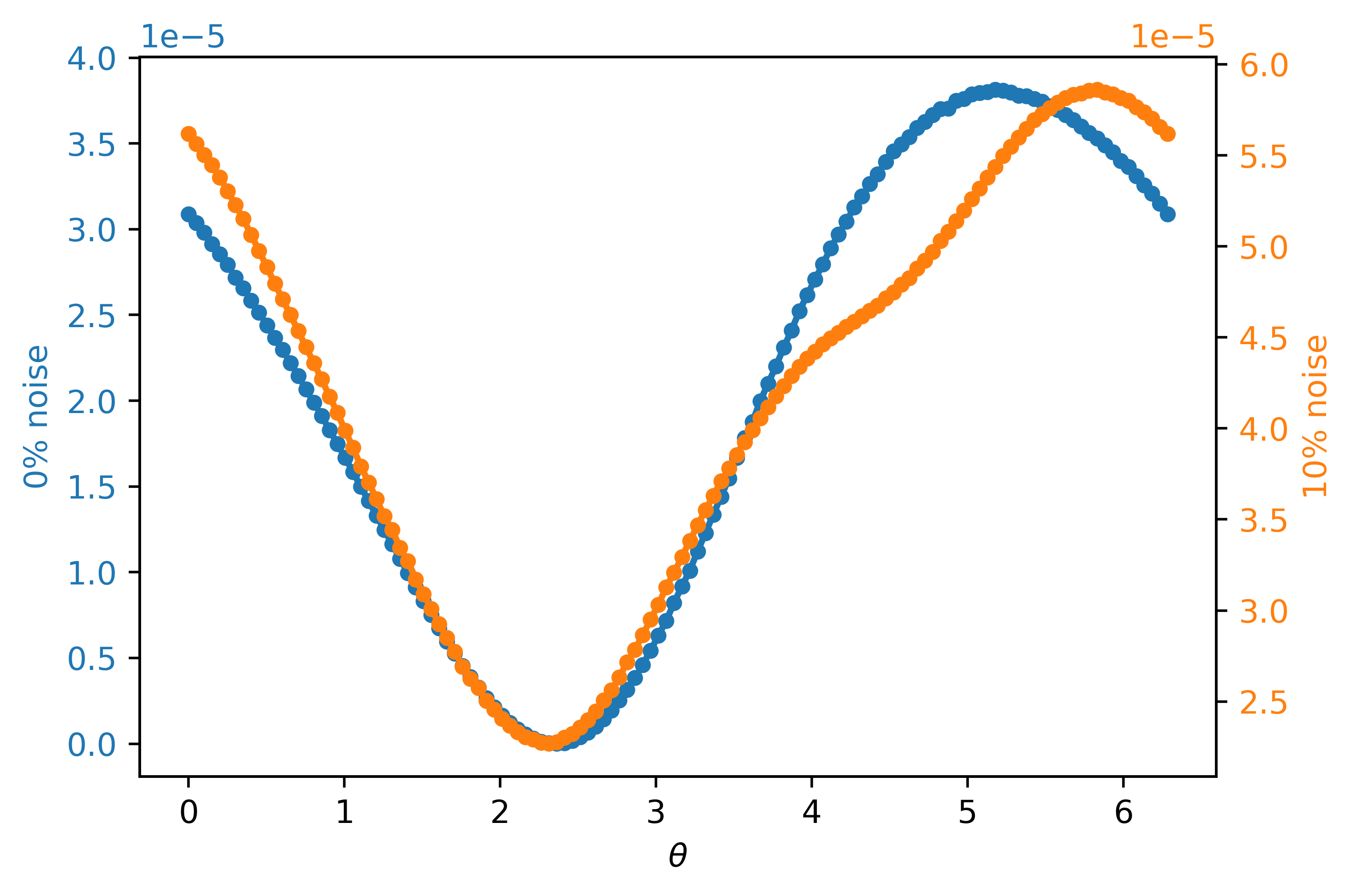}\label{W_comp}}
        \qquad
        \subfloat[$L^{2}$ landscape at R = 0.5]{\includegraphics[width=0.4\textwidth]{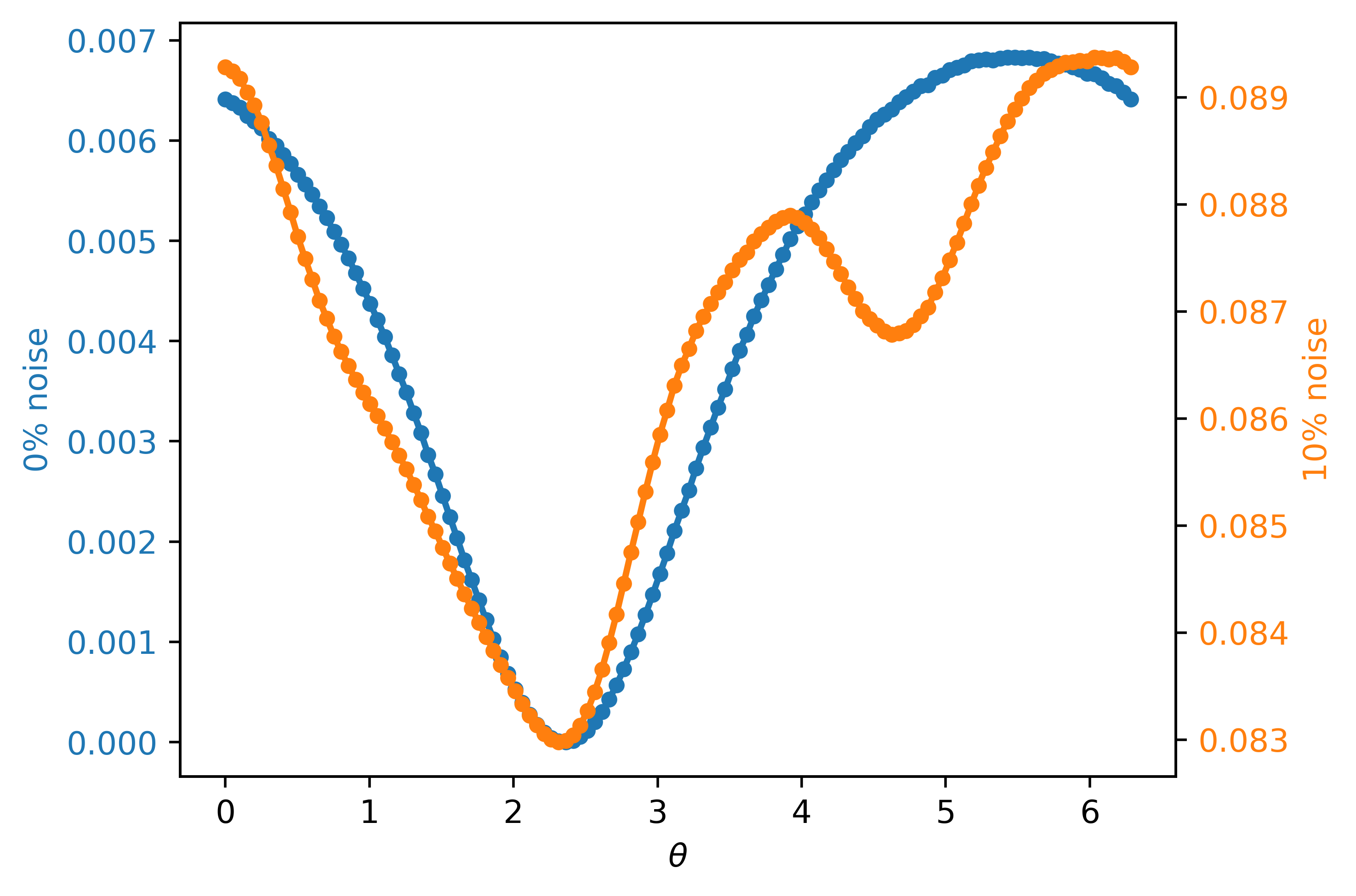}\label{L_comp}}
        \caption{The first row: landscapes of Example \ref{6_4} for $W_{2}$ distance (\ref{W_land}) and $L^{2}$ distance (\ref{L_land}) with $10\%$ noise in the data. The second row: cross sections of \ref{W_land} and \ref{L_land} at R = 0.5.}\label{fig5}
    \end{figure}
\end{example}

There are several points to note in the experiments. During the inversion, we take the first $m$ iterations using the $L^{2}$ distance without any regularization, where $m$ ranges from 5 to 15. During the initial iterations, only a general contour is captured, and the image has not been significantly interfered by the noise. Even if some noise information is fitted, it will not be severely overfitted in the next stage but will stay at an acceptable level due to the different noise performance of $W_{2}$ and $L^{2}$. Thus we consider using the $L^{2}$ distance for the preliminary optimization to save time. This procedure also allows the amplitude to reach a suitable value in short time. Not only that, it efficaciously prevents the occurrence of low-frequency artifacts under the $W_{2}$ distance. Another important issue is the choice of $a$ in (\ref{normalization}). Within a reasonable range, the value of $a$ does not have an apparent impact on the reconstruction results. For the numerical examples above, we can choose $a=2$, yielding quite accurate reconstructions. Furthermore, $a$ can also be determined in terms of the range of each NtD data, which improves the resolution of the reconstruction. It is shown from the numerical results that the resolution goes down as $a$ increases.

\section{Conclusion }
We have presented an efficient algorithm for solving the inverse conductivity problem with the quadratic Wasserstein distance. Our method is based on the simplified formulation of the optimal transportation problem on $\mathbb{S}^{1}$. A general form of the Fréchet gradient is derived through classical theory in OT. Then $W_{2}$ is coupled to the Barzilai-Borwein gradient method to implement the optimization. Our approach significantly enhances the quality of the reconstruction at the cost of a slight increase in computation. Numerical results on several examples demonstrate that our method performs better in removing background noise and identifying the singular shape than the traditional regularization method. The algorithm shows a strong resilience to noise and obtains reasonably accurate reconstructions in terms of location, magnitude, and shape. In addition, the proposed method can tackle the problem of overfitting to a certain extent, which makes the inversion process more stable. Important low-frequency information can be effectively distinguished. To the best of our knowledge, this is the first attempt to apply the $W_{2}$ metric to solve the Calderón problem (EIT). The success of our method shows that $W_{2}$ is an excellent choice for solving the severely ill-posed inverse problem. 

We point out some future directions of this research. The numerical algorithm in this paper is aimed at optimal transportation problems on closed smooth curves, which corresponds to inverse problems in 2D. In the case of three-dimensional inverse problems, the data domain is no longer to curves, and new methods must be developed to solve OT on surfaces in $\mathbb{R}^{3}$. Another interesting future direction is to develop solid theoretical analysis and advanced numerical techniques for the $W_{2}$ based inversion process. Despite the favorable properties of $W_{2}$, it is still unclear how these properties manifest themselves in most problems, especially for the convexity property. Finally, for the EIT problem, it is natural to consider the boundary data with zero mean, which may not be satisfied for other inverse problems. It would be interesting to consider more general data sets. In the more general setting, appropriate data normalization strategies should be designed to maintain the structure of the data and the metric. It is plausible to employ other OT-based distances, such as KR-norm \cite{metivier2016optimal} and Wasserstein–Fisher–Rao metric \cite{zhou2018wasserstein} to address this issue.\\

\textbf{Acknowledgments} The work was supported in part by National Natural Science Foundation of China (11621101; U21A20425) and a Key Laboratory of Zhejiang Province.


\begin{thebibliography}{50}




\bibitem{abraham2017tomographic}
I. Abraham, R. Abraham, M. Bergounioux, G. Carlier, Tomographic reconstruction from a few views: a multi-marginal optimal transport approach, Applied Mathematics \& Optimization 75.1 (2017): 55-73.

\bibitem{adler2021electrical}
A. Adler, D. Holder, Electrical Impedance Tomography: methods, history and applications, CRC Press, 2021.

\bibitem{ambrosio2013user}
L. Ambrosio, N. Gigli, A user’s guide to optimal transport, Modelling and optimisation of flows on networks. Springer, Berlin, Heidelberg, 2013. 1-155.

\bibitem{bao2015inverse}
G. Bao, P. Li, J. Lin, F. Triki, Inverse scattering problems with multi-frequencies, Inverse Problems 31.9 (2015): 093001.

\bibitem{bao2020numerical}
G. Bao, X. Ye, Y. Zang, H. Zhou, Numerical solution of inverse problems by weak adversarial networks, Inverse Problems 36.11 (2020): 115003.

\bibitem{Benamou2000}
JD. Benamou, Y. Brenier, A computational fluid mechanics solution
to the Monge-Kantorovich mass transfer problem, Numerische Mathematik 84.3 (2000): 375-393.

\bibitem{benamou2014numerical}
JD. Benamou, BD. Froese, AM. Oberman, Numerical solution of the optimal transportation problem using the Monge–Ampère equation, Journal of Computational Physics 260 (2014): 107-126.

\bibitem{borcea2002electrical}
L. Borcea, Electrical impedance tomography, Inverse Problems 2002; 18(6):R99–R136.

\bibitem{brenier1991polar}
Y. Brenier, Polar factorization and monotone rearrangement of vector valued functions, Comm. Pure Appl. Math. 44 (1991), 375–417.


\bibitem{caffarelli1996boundary}
LA. Caffarelli, Boundary regularity of maps with convex potentials--II, Annals of Mathematics 144.3 (1996): 453-496.

\bibitem{calderon2006inverse}
AP. Calderón, On an inverse boundary value problem, Computational \& Applied Mathematics 25 (2006): 133-138.

\bibitem{chen2018quadratic}
J. Chen, Y. Chen, H. Wu and D. Yang, The quadratic Wasserstein metric for earthquake location, Journal of Computational Physics 373 (2018): 188-209.

\bibitem{cheney1999electrical}
M. Cheney, D. Isaacson and JC. Newell, Electrical impedance tomography, SIAM review 41.1 (1999): 85-101.

\bibitem{cheney1990noser}
M. Cheney, D. Isaacson, JC. Newell, S. Simske, and J. Goble, NOSER: An algorithm for solving the inverse conductivity problem, International Journal of Imaging systems and technology 2.2 (1990): 66-75.


\bibitem{chung2005electrical}
ET. Chung, TF. Chan, XC. Tai, Electrical impedance tomography using level set representation and total variational regularization, Journal of Computational Physics 205.1 (2005): 357-372.

\bibitem{Cuturi2013Sinkhorn}
M. Cuturi, Sinkhorn distances: Lightspeed computation of optimal transport, Advances in neural information processing systems 26 (2013).

\bibitem{delon2010fast}
J. Delon, J. Salomon, A. Sobolevski, Fast transport optimization for Monge costs on the circle, SIAM Journal on Applied Mathematics 70.7 (2010): 2239-2258.


\bibitem{engquist2014application}
B. Engquist, BD. Froese, Application of the Wasserstein metric to seismic signals, arXiv preprint arXiv:1311.4581 (2013).

\bibitem{engquist2020quadratic}
B. Engquist, K. Ren, Y. Yang, The quadratic Wasserstein metric for inverse data matching, Inverse Problems 36.5 (2020): 055001.

\bibitem{engquist2020optimal}
B. Engquist, Y. Yang, Optimal transport based seismic inversion: Beyond cycle skipping, Communications on Pure and Applied Mathematics (2021).

\bibitem{fan2020solving}
Y. Fan, L. Ying, Solving electrical impedance tomography with deep learning, Journal of Computational Physics 404 (2020): 109119.

\bibitem{figalli2011optimal}
A. Figalli, C. Villani, Optimal transport and curvature, Nonlinear PDE’s and Applications. Springer, Berlin, Heidelberg, 2011. 171-217.

\bibitem{Glimm2003optical}
T. Glimm, V. Oliker, Optical design of single reflector systems and the Monge–Kantorovich mass transfer problem, Journal of Mathematical Sciences 117.3 (2003): 4096-4108.

\bibitem{heaton2022wasserstein}
H. Heaton, SW. Fung, AT. Lin, S. Osher, W. Yin, Wasserstein-based projections with applications to inverse problems, SIAM Journal on Mathematics of Data Science 4.2 (2022): 581-603.

\bibitem{haker2004optimal}
S. Haker, L. Zhu, A. Tannenbaum, et al, Optimal mass transport for registration and warping, International Journal of computer vision 60.3 (2004): 225-240.

\bibitem{isaacson1986distinguishability}
D. Isaacson, Distinguishability of conductivities by electric current computed tomography, IEEE Transactions on Medical Imaging 5.2 (1986): 91-95.

\bibitem{jin2012reconstruction}
B. Jin, T. Khan, P. Maass, A reconstruction algorithm for electrical impedance tomography based on sparsity regularization, International Journal for Numerical Methods in Engineering 89.3 (2012): 337-353.

\bibitem{kantorovich1960mathematical}
LV. Kantorovich, Mathematical methods of organizing and planning production, Management Science 6.4 (1960): 366-422.

\bibitem{knowles1998variational}
I. Knowles, A variational algorithm for electrical impedance tomography, Inverse Problems 14.6 (1998): 1513.

\bibitem{kohn1990numerical}
RV. Kohn, A. McKenney, Numerical implementation of a variational method for electrical impedance tomography, Inverse Problems 6.3 (1990): 3

\bibitem{logg2012automated}
A. Logg, KA. Mardal, G. Wells, Automated Solution of
Differential Equations by the Finite Element Method, Vol. 84. Springer Science \& Business Media, 2012.

\bibitem{mccann2001polar}
RJ. McCann, Polar factorization of maps on Riemannian manifolds, Geometric \& Functional Analysis GAFA 11.3 (2001): 589-608.

\bibitem{metivier2019graph}
L. Métivier, R. Brossier, Q. Merigot, É. Oudet, A graph space optimal transport distance as a generalization of $L_{p}$ distances: application to a seismic imaging inverse problem,  Inverse Problems 35.8 (2019): 085001.

\bibitem{metivier2016optimal}
L. Métivier, R. Brossier, Q. Merigot, É. Oudet, An optimal transport approach for seismic tomography: application to 3D full waveform inversion, Inverse Problems 32.11 (2016): 115008.

\bibitem{monge1781memoire}
G. Monge, Mémoire sur la théorie des déblais et des remblais, Mem. Math. Phys. Acad. Royale Sci. (1781): 666-704.


\bibitem{neuberger2009sobolev}
J. Neuberger, Sobolev Gradients and Differential Equations, Springer Science \& Business Media, 2009.

\bibitem{peyre2019computational}
G. Peyré, M. Cuturi, Computational optimal transport: With applications to data science, Foundations and Trends® in Machine Learning 11.5-6 (2019): 355-607.

\bibitem{peyre2018comparison}
R. Peyre, Comparison between $W_{2}$ distance and $\dot{H}^{-1}$ norm, and localization of Wasserstein distance, ESAIM: Control, Optimisation and Calculus of Variations 24.4 (2018): 1489-1501.

\bibitem{rabin2011transportation}
J. Rabin, J. Delon, Y. Gousseau, Transportation distances on the circle, Journal of Mathematical Imaging and Vision 41.1 (2011): 147-167.

\bibitem{rockafellar1970convex}
RT. Rockafellar, Convex Analysis, Princeton University Press, 1970.

\bibitem{rondi2001enhanced}
L. Rondi, F. Santosa, Enhanced electrical impedance tomography via the Mumford–Shah functional, ESAIM: Control, Optimisation and Calculus of Variations 6 (2001): 517-538.

\bibitem{santambrogio2015optimal}
F. Santambrogio, Optimal Transport for Applied Mathematicians, Birkäuser, NY 55.58-63 (2015): 94.

\bibitem{solomon2015convolutional}
J. Solomon, F. De Goes, G. Peyré, M. Cuturi et al, Convolutional Wasserstein distances: Efficient optimal transportation on geometric domains, ACM Transactions on Graphics (ToG) 34.4 (2015): 1-11.

\bibitem{uhlmann2009electrical}
G. Uhlmann, Electrical impedance tomography and Calderón's problem, Inverse Problems 25.12 (2009): 123011.

\bibitem{villani2021topics}
C. Villani, Topics in Optimal Transportation, Vol. 58. American Mathematical Soc., 2021.

\bibitem{Wang2004on}
XJ. Wang, On the design of a reflector antenna II, Calculus of Variations and Partial Differential Equations 20.3 (2004): 329-341.

\bibitem{wexler1985impedance}
A. Wexler, B. Fry, MR. Neuman, Impedance-computed tomography algorithm and system, Applied Optics 24.23 (1985): 3985-3992.

\bibitem{yang2018application}
Y. Yang, B. Engquist, J. Sun, BF. Hamfeldt, Application of optimal transport and the quadratic Wasserstein metric to full-waveform inversion, Geophysics 83.1 (2018): R43-R62.

\bibitem{zhou2018wasserstein}
DT. Zhou, J. Chen, H. Wu, DH. Yang, LY. Qiu, The Wasserstein-Fisher-Rao metric for waveform based earthquake location, arXiv preprint arXiv:1812.00304 (2018).


\end{thebibliography}
\end{document}